\documentclass[11pt]{amsart}

\usepackage{graphics}
\usepackage{amscd}
\usepackage{amssymb}

\newcommand{\bz}{\mathbb{Z}}
\newcommand{\br}{\mathbb{R}}
\newcommand{\bc}{\mathbb{C}} 
\newcommand{\tp}{\widetilde{\bz_p}}

\newcommand{\fred}{S^2\times S^2}
\newcommand{\lan}{\left\langle}
\newcommand{\ran}{\right\rangle}
\newcommand{\tet}{\mbox{Tet}}
\newcommand{\oct}{\mbox{Oct}}
\newcommand{\icos}{\mbox{Icos}}
\newcommand{\aut}{\operatorname{Aut}}
\newcommand{\out}{\operatorname{Out}}
\newcommand{\inn}{\operatorname{Inn}}
\newcommand{\image}{\operatorname{im}}
\newcommand{\rk}{\operatorname{rk}}

\newcommand{\gkk}{g_{k^2}}
\newcommand{\gqq}{g_{q^2}}
\newcommand{\gqaqa}{g_{q_a^2}}
\newcommand{\gqbqb}{g_{q_b^2}}
\newcommand{\matrixa}{\left(\begin{smallmatrix} -1&0\\
	0&-1\end{smallmatrix}\right)}
\newcommand{\matrixb}{\left(\begin{smallmatrix} 0&1\\
	1&0\end{smallmatrix}\right)}

\newenvironment{eg}%
  {
    \begin{trivlist}{}{\setlength{\leftmargin}{0pt}}%
    \item{\textbf{Example}\enspace}\par%
  }%
  {
    \end{trivlist}%
  }
\newenvironment{defn}%
  {
    \begin{trivlist}{}{\setlength{\leftmargin}{0pt}}%
    \item{\textbf{Definition}\enspace}%
  }%
  {
    \end{trivlist}%
  }

\newenvironment{revision}{}{}

\title[Groups which act pseudofreely on $S^2\times S^2$. ]{Groups which act pseudofreely on $S^2\times S^2$. }

\author[McCooey]{Michael P. McCooey}

\address{Department of Mathematics, Franklin and Marshall College\\
P. O. Box 3003\\
Lancaster, PA , 17604-3003 } 
\email{mmccooey@member.ams.org}

\urladdr{http://edisk.fandm.edu/michael.mccooey/}

\newtheorem{thm}{Theorem}[section]
\newtheorem{prop}[thm]{Proposition}
\newtheorem{lem}[thm]{Lemma}
\newtheorem{cor}[thm]{Corollary}
\theoremstyle{definition}

\theoremstyle{remark}

\begin{document}
\begin{abstract}
Recall that a {\em pseudofree} group action on a space $X$ is one whose 
set of singular orbits forms a discrete subset of its orbit space.
Equivalently  -- when $G$ is finite and $X$ is compact -- 
the set of singular points in $X$ is finite. In this paper,
we classify all of the finite groups which admit pseudofree
actions on $\fred$. The groups are exactly those which admit
orthogonal pseudofree actions on $\fred\subset\br^3\times\br^3$,
and they are explicitly listed.

This paper can be viewed as a companion to a preprint of 
Edmonds,
which uniformly treats the case in which the second Betti number of 
a four-manifold $M$ is at least three.
\end{abstract}
Ê
\maketitle
Ê\section{Introduction}


	The classification of free actions of groups on spheres
and products of spheres is an important problem
in topology. In even dimensions, however, the question of which
\emph{groups} admit free actions is much less interesting: It follows
from the Lefschetz fixed point theorem that such a group must
admit a faithful, and rather special, representation on 
homology. Thus, for example, $\bz_2$ is the only group which can 
act freely on $S^{2n}$, or freely and
orientably on $S^{2n}\times S^{2n}$. 
	The Lefschetz fixed point theorem  puts the same sort
of restriction on free actions of groups on closed, simply-connected
four-manifolds: in the homologically trivial case, there are none.
This is one motivation for the concept of  \emph{pseudofree}
actions. 

An action of a finite group $G$ on a space $X$ is pseudofree
if it is free on the complement of a discrete set of points.\footnote{Some 
authors require that a pseudofree action also
be \emph{semifree}, so that each singular point is fixed by the
entire group. 
}
Such actions make the singular set of the group action as small
as is compatible with the Lefschetz theorem. It is a theorem of
Edmonds~\cite{Actions} that if $b_2(M^4)\ge 3$, 
the only groups which
can act pseudofreely on $M$ and trivially on $H_*(M)$ are cyclic. 
In this paper, we treat the case $M =\fred$ -- which is especially interesting in its own right -- for all possible actions on homology.

Let $G$ be a finite group with an action on 
$S^2\times S^2$ which 
preserves orientation and is locally linear and pseudofree.
Two questions arise naturally:   Which
groups $G$ admit such actions? And given a group which acts pseudofreely, can we
classify the actions?  A partial answer to the latter
question is known: if $G$ is cyclic,
the invariants of Wilczy\'nski and Bauer~\cite{B&W}, \cite {Wil}, 
detect whether two given pseudofree actions (on any
simply-connected four-manifold)  are topologically conjugate. 
The work of Edmonds and Ewing~\cite{EdmondsEwing}, as well
as that of other authors, addresses the
question of which combinations of these invariants are
realizable, but there are number-theoretical difficulties involved in
making the answer explicit.
	
Evidence  from the study of other simple $4$-manifolds 
(See~\cite{CS}, \cite{Kulkarni} for the case of $S^4$, \cite{Wil2}, 
\cite{HL} for $\bc P^2$) hints at an answer to the first question:
the groups which admit pseudofree actions are exactly those
which can also act pseudofreely and orthogonally on 
$S^2\times S^2\subset \br^6$. Our main result is that this is
indeed the case: we classify the groups acting orthogonally
and pseudofreely, and prove that these are the only groups
admitting even locally linear pseudofree actions. (More precisely,
the groups and their corresponding representations on homology
are exactly those which occur in the linear case.)
The remaining problem, that of classifying the actions themselves,
is left for future work.

Here is a brief outline of the paper:
The homology representation of a group acting on $\fred$ gives rise
to a short exact sequence.  Certain group extension problems 
related to this sequence are solved, and then the geometry
of linear actions is used to determine which of these extensions 
admit pseudofree orthogonal actions.
Actions are constructed in the course of
the proofs. The classification of groups acting
linearly and pseudofreely is stated explicitly in 
Theorem~\ref{biglist}, and is reformulated in terms of the homology
representation in Corollary~\ref{linearclassification}.
In Section~\ref{nonlinearcase}, the conditions of 
Corollary~\ref{linearclassification} are shown to be 
necessary in the locally linear case, as well, thereby
proving the main theorem. The proof uses a series of lemmas which
apply orbit-counting arguments to the singular set of a 
group action, and then some group-theoretical and cohomological 
calculations to rule out minimal potential pathologies.

\section{Group theory}\label{grouptheory}
Suppose a group $G$ acts on $\fred$. The induced 
action of $G$ on $H_2(S^2\times S^2)$ defines a 
representation $\varphi : G \rightarrow GL(2, \bz)$. Since 
$\varphi$ must respect the intersection form, it must leave the 
the positive and negative definite subspaces of $H_2(\fred)$
invariant. $H_{2}^{+}$ is spanned (rationally) by
$x+y$, and $H_{2}^{-}$ is spanned by $x-y$, where $x$ and $y$
represent the standard generators. With respect to  this basis,
$\varphi(g)$ must have the form 
$\left( \begin{array}{cc}
	\pm 1 & 0\\
	0 & \pm 1
	\end{array}\right)$ for any $g\in G$.
It follows that with  respect to the standard basis, 

$$\varphi(g) \in \left\langle
	\pm\left(\begin{array}{cc}1&0\\0&1\end{array}\right),
	\pm\left(\begin{array}{cc}0&1\\1&0\end{array}\right)\right\rangle
	\cong \bz_2\times\bz_2.$$

Thus the representation on homology induces a short exact sequence

\begin{equation*}\label{sequence}
1\rightarrow K\rightarrow G\stackrel{\varphi}{\rightarrow}
Q\rightarrow 1,
\end{equation*}
where $K$ acts trivially on homology and $Q\subset \bz_2\times \bz_2$.
This sequence allows us to approach general $G$-actions by first 
considering homologically trivial actions and then dealing with
associated extension problems. 
In this section, we address the group theory involved
in the extension problems,
and also some cohomology calculations which will be needed later on
in the nonlinear case. 
The results, though technical, are for the most part routine.

First recall some generalities about the classification of group extensions:

If $K$ is an abelian group, and $1\rightarrow K\rightarrow G
\rightarrow Q\rightarrow 1$ is an extension, then $G$ acts on $K$ by
conjugation. Since $K$ acts trivially on itself, conjugation induces
a well-defined action of $Q$ on $K$.
 For a specific $g\in G$, the conjugation automorphism
$h\mapsto g^{-1}hg$  will be denoted $\mu_g$. The action of $Q$
defines a homomorphism  $\psi:Q\rightarrow \aut(K)$ and a 
$Q$-module structure
on $K$. As is well known, the extensions $1\rightarrow K\rightarrow 
G\rightarrow Q\rightarrow 1$ which give rise to the action $\psi$ 
of $Q$ on $K$ are classified by $H^2(Q; K)$.

If $K$ is nonabelian, the situation is a bit more delicate:
 an extension problem 
$1\rightarrow K\rightarrow G\rightarrow Q\rightarrow 1$ with
nonabelian kernel is described by the ``abstract kernel''
$(K, Q, \psi)$, where $\psi: Q\rightarrow \out(K)$ describes the 
``outer action'' of $Q$ on $K$. For general $(K, Q, \psi)$, an
extension might not exist: $\psi$ determines an obstruction
cocycle in $H^3(Q; Z(K))$. The obstruction measures, roughly, 
whether it is possible to simultaneously realize the (outer)
actions of each element of $Q$ by conjugations in $G$. If the
obstruction vanishes, then the extensions realizing $(K, Q, \psi)$
are in (non-canonical) 1-1 correspondence with $H^2(Q; Z(K))$.
For detailed accounts of this theory, see \cite{Brown} or \cite{MacLane}.

In this section, we will  mainly be concerned with extensions
with quotient $\bz_2$. Assume such an extension exists, and
pick some $q\in G\setminus K$.
Two pieces of data determine the structure of $G$:
\begin{enumerate}
\item{The particular value in $K$ of the square of $q$. Denote it
by $\gqq$. }
\item{The automorphism $\mu_q$ of $K$. In the abstract (when a 
specific extension is not given in advance), this automorphism
will be denoted $\gamma$.}
\end{enumerate}
The choices of $\gqq$ and $\gamma$ are not arbitrary: 
$\gamma$ must always fix $\gqq$, and $\gamma^2$ must be
an inner automorphism of $K$. In specific instances, there may 
be other restrictions as well.

\goodbreak

\noindent {\bf Notation:}  The dihedral group of order $2m$ will be
denoted $D_m$. Similarly,  the binary dihedral group of 
order $4m$, with presentation $\lan k, q|k^m= q^2, q^{-1}kq=k^{-1}\ran$,
will be denoted $D^*_m$.
(This group  is also commonly referred to as a dicyclic or  generalized
quaternion group.)

\begin{lem}\label{Hcycliclemma}

Let $K\cong \bz_n$.  Each  extension 
$0\rightarrow K\rightarrow G\rightarrow \bz_2\rightarrow 1$
is of one of the following forms:  
\begin{enumerate}
\item{
	$\bz_n\rtimes\bz_2$, where the semidirect product automorphism
is given by a certain tuple $(\delta_2, \epsilon_2, \ldots, 
\epsilon_{p_i},\ldots, \epsilon_{p_k})$. Every split extension is
of this type.}

\item{
	$(\bz_{m_-}\rtimes\bz_{2^{n_2+1}})\times \bz_{m_+}$, where the
semidirect product automorphism is inversion.} 
\item{
	$D^*_{2^{n_2-1}m_-}\times \bz_{m_+}$.}

\end{enumerate}
If $n$ is odd, then every extension is split. Notation is explained 
in the proof.

\end{lem}

\begin{proof}
We use additive notation for the group operation in $\bz_n$.

Let $n=p_1^{n_{p_1}}\cdots p_q^{n_{p_k}}$ be a prime 
factorization. Then
$\bz_n\cong\bz_{p_1^{n_{p_1}}}\times\ldots\times\bz_{p_q^{n_{p_k}}}$, and
$\gamma$ restricts to an automorphism of each of the factors.

If $p$ is odd, $\aut(\bz_{p^{n_p}})$ is cyclic, so each order two
automorphism of $\bz_{p^{n_p}}$ sends a generator to its inverse. 
However, $\aut(\bz_{2^{n_2}})\cong\bz_2\times\bz_{2^{n_2-2}}$, so
if $n_2>2$, $\bz_{2^{n_2}}$ has three automorphisms of order two:
$x\mapsto -x$, $x\mapsto(1+2^{n_2-1}) x$, and $x\mapsto
(-1+2^{n_2-1})x$. Thus $\gamma |_{\bz_{2^{n_2}}}$ is encoded by
a pair $(\delta_2, \epsilon_2)$, where $\delta_2=1 \mbox{ or } 0$, 
$\epsilon_2=\pm 1$, and  $(\delta_2, \epsilon_2):
x\mapsto (\epsilon_2+\delta_2\cdot 2^{n_2-1})x$. For $p>2$,
we simply have $\gamma|_{\bz_{p^{n_p}}}:x\mapsto \epsilon_p\cdot x$.
With this notation, we write $\gamma =(\delta_2, \epsilon_2, \ldots, 
\epsilon_{p_i},\ldots, \epsilon_{p_k})$. 
Then $\bz_n$ can be viewed as a product 
	$$\bz_n\cong\bz_{2^{n_2}}\times\bz_{m_-}\times\bz_{m_+},$$
where $\bz_{m_-}$ is the (odd) ``$-1$-eigenspace'' for $\gamma$,
and $\bz_{m_+}$ is the (odd) ``$+1$-eigenspace''.

Since $\gamma$ fixes $\gqq$, $\gqq\equiv 0 \pmod{p^{n_p}}$ 
for each odd 
$p$ with $\epsilon_p=-1$, so $\gqq\equiv 0 \pmod {m_-}$.
On the other hand, any element of $\bz_{m_+}$ is a  multiple of 2, 
so $q$ may be normalized to be trivial mod $m_+$, as well.
We may therefore assume that $\gqq\in \bz_{2^{n_2}}\times 0\times 0$,
and we have a subextension $1\rightarrow \bz_{2^{n_2}}\rightarrow G'
\rightarrow \bz_2\rightarrow 1$. These are classified by $H^2(\bz_2;
\bz_{2^{n_2}})$, where $(\delta_2, \epsilon_2)$ defines the module
structure on $\bz_{2^{n_2}}$. Brown~\cite[IV.4.2]{Brown} computes
that if $\delta_2=1$, then $H^2(\bz_2; \bz_{2^{n_2}})=0$, so
every extension is split. Thus if $\delta_2=1$, the main extension
$1\rightarrow\bz_n\rightarrow G\rightarrow\bz_2\rightarrow 1$ is a
semidirect product. 
Henceforth we assume $\delta_2=0$.
We now have $H^2(\bz_2; \bz_{2^{n_2}})\cong \bz_2$, so
$H^2(\bz_2; \bz_n)\cong \bz_2$, as well. The split extensions are
again semidirect products. But for each $\gamma$, there is one 
non-split extension.
 Since $\bz_{m_+}$ is central in $G$, the 
sequence can be written
$$0\rightarrow (\bz_{2^{n_2}}\times\bz_{m_-})\times\bz_{m_+}
\rightarrow G''\times\bz_{m_+}\rightarrow\bz_2\rightarrow 1,$$
so we may pretend for the moment that $\bz_{m_+}=\{0\}$, and focus
attention on $G''$.  

The structure of $G''$ depends on the sign of $\epsilon_2$. If 
$\epsilon_2=1$, then $\bz_{2^{n_2}}$ is central, and $G''\cong
\bz_{m_-}\rtimes\bz_{2^{n_2+1}}$. If $\epsilon_2=-1$, then 
$G''\cong D^*_{2^{n_2-1}m_-}$. By appending the $\bz_{m_+}$ factor,
we recover $G$ in each case.

\end{proof}

Next we consider the possibility that $K$ is a dihedral group. 
Since extensions of dihedral groups have nonabelian kernel,
the algebra is a little more technical than in the cyclic case.

Given an automorphism $\gamma$ of $D_n$ whose square
is inner, there exists $\gqq\in D_n$ such that $\gamma^2=\mu_{\gqq}$, 
and this choice of $\gqq$ is unique up to a factor in $Z(D_n)$.
$Z(D_n)$ is trivial if $n$ is odd, and has order 2 if $n$ is even.
Clearly, $\gqq$ is fixed by $\gamma^2$.
Working through the definitions in~\cite{Brown} or \cite{MacLane} 
shows that the obstruction to realizing the abstract kernel 
$(D_n, \bz_2, \gamma)$ as an extension  vanishes exactly if $\gqq$ 
is fixed by $\gamma$. It is straightforward to verify that vanishing 
of the obstruction depends only on the outer automorphism class 
of $\gamma$.

\begin{revision}

\noindent {\bf Claim 1:} $\aut(D_n)\approx\begin{cases}S_3 & \text{if } n=2,\\
 \bz_n\rtimes \aut(\bz_n)& \text{if } n>2\end{cases}$
\end{revision}

We use the presentation $D_n\approx\lan s,t|s^n=t^2=1,
tst=s^{-1}\ran$.

\begin{revision}
If $n=2$, then explicitly checking that the 2-cycle $(s,t)$ and the 3-cycle $(s,t,st)$ respect the group operation shows that any permutation of $\{s,t,st\}$ defines an automorphism of $D_2$.

If $n>2$ then, since an automorphism must send $s$ to another element of order $n$,  
and $t$ to an element which does not commute with $f(s)$, \end{revision}
any automorphism is of the form $f_{a,b}$, where
$f_{a,b}(t)=s^at$, and $f_{a,b}(s)=s^b$, for $a\in \bz_n$ and 
$(b,n)=1$. By calculating $f_{a,b}\circ f_{c,d}=f_{bc+a, bd}$, we
see that $\aut(D_n)$ is a semidirect product, as claimed.

\begin{revision}Let us first deal with the case $n=2$.
Since $D_2$ is abelian, $\gamma^2=1$.

\noindent{\bf Case 1:} $\gamma$ is nontrivial. 

For convenience, we choose generators $a, b$ of $D_2$,  and assume that $\gamma$
transposes them, leaving $ab$ fixed. Since $\gqq$ must be fixed by $\gamma$, $\gqq=ab$ or $1$; replacing $g_q$ by $aq$ if necessary, we may assume $q^2=1$. Then $G=\lan a, b, q\ |\ a^2=b^2=(ab)^2=q^2=1,qaq^{-1}=b \ran$, 
or more simply, setting $s=aq, t=b$, $G=\lan s, t\ |\ s^4=t^2=1, tst=s^{-1}\ran$, with
$D_2$ included as the subgroup $\lan s^2, t\ran$.

\noindent {\bf Case 2:} $\gamma$ is trivial.

In this case, $G$ is abelian. It is straightforward to check that $G=\bz_4\times\bz_2$ or $G=\bz_2\times\bz_2\times\bz_2$.

Assume henceforth that $n>2$.

\end{revision}

\noindent {\bf Claim 2:} $\inn(D_n)\approx 2\bz_n\rtimes\{\pm 1\}\subset
\bz_n\rtimes \aut(\bz_n)$. 

It suffices to check that $\mu_{s^a}=f_{2a, 1}$, and $\mu_t=f_{0, -1}$.

Now let $\gamma=f_{a,b}$. Then $\gamma^2=f_{ab+a, b^2}.$ If
$\gamma^2\in\inn (D_n)$, then $b^2=\pm 1$. But if $b^2=\pm 1$, then 
$b+1$ must be a multiple of 2 in $\bz_n$, so $ab+a\in 2\bz_n$ for 
any $a\in\bz_n$. Thus $\gamma^2$ is inner iff $b^2=\pm 1$. Since
$\gamma^2$ must be inner for an extension to exist, we assume
henceforth that $b^2=\pm 1$.

\noindent{\bf Claim 3:} 
If $n$ is odd, the obstruction always vanishes.
If $n$ is even and $b^2\equiv 1\ (n)$, the obstruction
vanishes unless $a$ and $\frac{b^2-1}{n}$ are both odd in $\bz$.
If $n$ is even and $b^2\equiv -1\ (n)$, the obstruction
vanishes unless $a$ and $\frac{b^2+1}{n}$ are both odd in $\bz$.

Suppose $n$ is even. Then the equation $2b'=b+1$ has two solutions 
mod $n$. Fix one. The other is $b'+\frac{n}{2}$.

\noindent{\bf Case 1:} $b^2=1$.

In this case, $\gamma^2=\mu_{(s^{ab'})}$, so 
$\gqq=s^{ab'}$ or $s^{a(b'+n/2)}$. If the first is fixed by
$\gamma$, the second is, as well. So we ask: when is
$\gamma(s^{ab'})=s^{ab'}$?

Since $\gamma(s^{ab'})=s^{bab'}$, $\gqq$ is fixed by $\gamma$
iff $(b-1)ab'\equiv 0\ (n)$. Notice that $2(b-1)ab'=a(b^2-1)
\equiv 0\ (n)$. However, it is possible that 
$(b-1)ab'\equiv \frac{n}{2}\ (n)$.
This occurs if $2(b-1)ab'=a(b^2-1)$ is an odd multiple of
$n$, which occurs only if $a$ and $\frac{b^2-1}{n}$ are both odd. 

\noindent{\bf Case 2:} $b^2=-1$.

Now $\gamma^2=\mu_{(s^{ab'}t)}$, and $\gqq=s^{ab'}t$ or
$s^{a(b'+n/2)}t$. $\gamma$ fixes $\gqq=s^{ab'}t$ iff $ab'\equiv
bab'+a\ (n)$, iff $(b-1)b'a+a\equiv 0\ (n)$. 
As above, everything works unless $2((b-1)b'a +a)=(b^2-1)a+2a
=(b^2+1)a$ is an odd multiple of $n$. This occurs only if
$a$ and $\frac{b^2+1}{n}$ are both odd.

If $n$ is odd, then similar but easier calculations show that
$\gqq$ is always fixed by $\gamma$. This establishes the claim.

\begin{defn} Let us say that an automorphism $\gamma$ of $D_n$ is 
\emph{admissible} if an extension 
$1 \rightarrow D_n \rightarrow G\rightarrow \bz_2\rightarrow 1$
exists with $\mu_k=\gamma$ for some  $k\in G\setminus D_n$. 
\end{defn}

We have just seen that \begin{revision} for $n>2$, \end{revision}$\gamma = f_{a, b}$ is admissible if and
only if $b^2\equiv\pm 1\ (n)$ and the obstruction mentioned in
claim~3 vanishes.

When extensions exist, they are in correspondence with
$$H^2(\bz_2, Z(D_n))\approx \left\{	
	\begin{array}{ll}
	0 & \mbox{ if $n$ is odd,}\\
	\bz_2 & \mbox{ if $n>2$ is even}
	\end{array}\right. $$

If $n$ is even, the extensions correspond to the two choices
of $\gqq$, which differ by the nontrivial element in $Z(D_n)$.
Thus
\begin{eqnarray*}
G &\cong & \lan s,t,q | s^n=t^2=1, tst=s^{-1}, q^2=\gqq, 
		q^{-1}sq=\gamma(s), q^{-1}tq=\gamma(t)\ran \\
  &\cong &\lan D_n, q | q^2=\gqq, q^{-1}gq=\gamma(g) \mbox{ for }
		g\in D_n\ran .
\end{eqnarray*}

This presentation depends on a particular $\gamma\in [\gamma]\in
\out(D_n)$. For some purposes, we might wish to normalize $\gqq$.
To this end, note that for $g\in D_n$, 
$(\mu_g\circ\gamma)^2=\mu_{g\gamma(g)}\circ\gamma^2 =
\mu_{g\gamma(g)\gqq}$. Thus $\gqq$ can be modified by left 
multiplication by any element of the form $g\gamma(g)$.
The sequence will split when one of these candidates for
$\gqq$ is the identity. 

To summarize:
\begin{lem}\label{dihedralextensionlemma}
Let $\gamma\in\aut (D_n)$. If $\gamma$ is admissible, there are extensions
$$1\rightarrow D_n\rightarrow G\rightarrow \bz_2\rightarrow 1$$
in which the conjugation action of some $q\in G\setminus D_n$
is given by $\gamma$.  
\begin{revision} If $n=2$, the extensions are $$1\rightarrow D_2\rightarrow D_4\rightarrow \bz_2\rightarrow 1,$$ with $D_2$ mapping to $\lan s^2, t\ran\subset D_4$, and the two abelian extensions with $G\approx \bz_4\times\bz_2$ and $G\approx(\bz_2)^3$. For $n>2$,
\end{revision}there is exactly one such extension
if $n$ is odd, and two if $n$ is even, and every extension is of this
form. The resulting groups $G$ have presentations of the form 
$$\lan D_n, q | q^2=\gqq, q^{-1}gq=\gamma(g) \mbox{ for } g\in D_n\ran.$$

\end{lem}

Let $\tet$, $\oct$, and $\icos$ denote the three symmetry groups of
the Platonic solids. 

\begin{lem}
\label{Platosplitslemma}
Let $K$ be one of the groups $\tet$, $\oct$, or $\icos$. Then every extension
$$1\rightarrow K \rightarrow G \rightarrow \bz_2\rightarrow 0$$ is split.
\end{lem}

\begin{proof}
Recall that $\tet\approx A_4$, $\oct\approx S_4$, and $\icos\approx A_5$. 
A group $G$ is said to be \emph{complete} if $Z(G)$ is trivial and every 
automorphism of 
$G$ is inner. For $n \ne 2,6$, $S_n$ is complete (see~\cite[th. 7.5]{Rotman}).
 Alternating groups are never
complete: for $n>2$, $A_n\subseteq S_n$, and conjugation by the transposition
$(12)$ is an automorphism of $A_n$ which is not inner in $A_n$.

\noindent{\bf Claim:} $\aut(A_4)\approx S_4$.

We know $S_4\subseteq \aut(A_4)$, since any conjugation by an element of
$S_4$ leaves $A_4$ invariant, and since $Z(S_4)$ is trivial. $A_4$ has 
the presentation $\lan x,y|x^3=y^2=(xy)^3=1\ran$. Any automorphism 
must send $x$ to an element of order 3, and $y$ to an element of 
order 2. Since $A_4$ contains 8 elements of order 3 and 3 elements of 
order 2, we see that $|\aut(A_4)|\le 3\cdot 8 =24$. But $|S_4|=24$.

\noindent{\bf Claim:} $\aut(A_5)\approx S_5$.

This follows by an argument similar to the preceding one: the 
presentation $\lan x,y|x^2=y^5=(xy)^3=1\ran$ gives an upper bound
$|\aut(A_5)|\le 180$. But $\aut(A_5)\supseteq S_5$, a group of
order 120.

Thus $\out(\tet)\approx\out(\icos)\approx \bz_2$, with the outer 
automorphism realized by conjugation by a transposition in $S_4$ or 
$S_5$. $\out(\oct)$ is trivial, since $S_4$ is complete.

Now, each of $\tet$, $\oct$, and $\icos$ has trivial center, so
in each case, $H^2(\bz_2; Z(K))=0$. Thus if an abstract kernel
admits a realization, it will be unique. Extensions do exist in
all cases: the two extensions of $\tet$ have $G\cong S_4$ and 
$G\cong \tet\times\bz_2$; the extensions of $\icos$
have $G\cong S_5$ and $G\cong \icos\times\bz_2$; and the unique
extension of $\oct$ has $G\cong\oct\times\bz_2$. All are
split extensions.

\end{proof}

The last classification lemma is somewhat more technical than
the others, but it turns out to be exactly what is necessary
when $\varphi(G)=\bz_2\times\bz_2$. 

\begin{lem}\label{technical}
Suppose $\lan a,b\ran\cong\bz_2\times\bz_2$, and $1\rightarrow \bz_n
\rightarrow G\stackrel{\varphi}{\rightarrow}\lan a, b\ran\rightarrow 1$ is
exact. Let $G_a=\varphi^{-1}(\lan a\ran )$ and 
$G_b=\varphi^{-1}(\lan b\ran )$. Finally, let $\gamma=\mu_{q_b}$, viewed
as an automorphism of $G_a$. Using the notation of 
Lemma~\ref{Hcycliclemma}, we have $\gamma |_{\bz_n} =(\delta_2, \epsilon_2, 
\ldots, \epsilon_{p_i},\ldots, \epsilon_{p_k})$.
If $G_a$ is abelian and
$1\rightarrow G_a\rightarrow G\rightarrow \lan b\ran\rightarrow 1$
does not split, then one of the following holds:
\begin{enumerate}
\item{$G\cong D^*_{2^{n_2}m_-}\times \bz_{m_+}$, with $\epsilon_2=-1$,
$G_a\cong\bz_{2n}$, and $\gamma(q_a)=q_a^{-1}.$}
\item{$G\cong((\bz_{m_-}\rtimes\bz_{2^{n_2+1}})\times \bz_{m_+})\rtimes
\bz_2$, with $G_a\cong\bz_a\times\bz_2$, $n_2>1$, $\epsilon_2=1$, and
$\gamma(q_a)=((n/2)k)q_a$.}
\item{$G\cong(D^*_{2^{n_2-1}m_-}\times\bz_{m_+})\times\bz_2$, with 
$G_a\cong \bz_n\times\bz_2$, $\epsilon_2=-1$, and $\gamma(q_a)=q_a$.}
\item{$G\cong(\bz_{m_-}\rtimes\bz_{2^{n_2+1}})\times\bz_{m_+}\times\bz_2$,
with  $G_a\cong \bz_n\times\bz_2$, $\epsilon_2=1$, and $\gamma(q_a)=q_a$.}
\end{enumerate}
In the latter two cases, $q_a$ and $q_b$ commute. Notation,
including some important normalizations for $q_a$ and $q_b$,
 is explained in the proof.
\end{lem}

\begin{proof} We use additive notation for the group operation in
$K$, but multiplicative notation in the (possibly nonabelian) group $G$.
Let $k$ generate $K=\bz_n$. $G$ is generated by $k$, together with
$q_a$ and $q_b$, where $\varphi(q_a)=a$ and $\varphi(q_b)=b$.
We assume $q_a$ and $q_b$ are normalized so that $\gqaqa$ and
$\gqbqb$ lie  in the Sylow 2-subgroup of $K$, and 
\begin{enumerate}
\item{$\gqaqa$ is either 0, or a generator of this subgroup.}
\item{If $\epsilon_2=1$, then $\gqbqb$ is either 0, or a generator.}
\item{If $\epsilon_2=-1$, then either $\gqbqb=0$, or it has order two.}
\end{enumerate}
With this in mind, the problem easily reduces to the case $n=2^{n_2}$,
since the restriction of $\gamma$ to the subgroup $\bz_{m_-m_+}$ of
$\bz_n$ is known. We assume henceforth that $n$ is a power of 2.

If $\gqbqb=0$, the sequence splits, so we assume $\gqbqb\ne 0$. Once
$G_a$ is known, $G$ is determined by the automorphism $\gamma$ of
$G_a$. In the proof of Lemma~\ref{Hcycliclemma}, we observed that
$\gamma|_{\bz_n}$ is described by a pair $(\delta_2, \epsilon_2)\in
\{0,1\}\times\{\pm 1\}$, where $k\mapsto (\epsilon_2 +\delta_2\cdot
2^{n_2-1})k$, and also that if $\delta_2=1$, then $1\rightarrow\bz_n
\rightarrow G_b\rightarrow \bz_2\rightarrow 1$ splits. So we may assume
$\gamma(k)=\epsilon_2\cdot k$ and consider the various cases for 
$\gamma(q_a)$. Note that $\gamma(q_a)=cq_a$ for some $c\in K$, and 
that $\gamma(q_a^2)=2c +\gqaqa$. Hence $c$ is determined mod $2^{n_2-1}$
by $\gamma|_{\bz_n}$. 

If $G_a\cong\bz_{2n}$, then the argument cited in Lemma~\ref{Hcycliclemma}
shows that $1\rightarrow G_a\rightarrow G\rightarrow \lan b\ran\rightarrow 1$
splits if $c\ne 0$. So we assume $\gamma(q_a)=q_a^{\pm 1}$, depending on 
$\epsilon_2$.  If $\epsilon_2=-1$, then $G\cong D^*_{2^{n_2}}$. 
If $\epsilon_2=1$, then $G\cong\bz_{2n}\times\bz_2$, and $(q_aq_b^{-1})^2=1$,
so $1\rightarrow G_a\rightarrow G\rightarrow \lan b\ran\rightarrow 1$ splits.

If $G_a\cong \bz_n\times\bz_2$, we have several cases:
\begin{enumerate}
\item{ If $c=0$ and $\epsilon_2=1$, then $G\cong\bz_{2^{n_2+1}}\times\bz_2$.
$G$ contains three elements of order 2, but all are contained in $G_a$,
so the sequence $1\rightarrow G_a\rightarrow G\rightarrow\lan b\ran
\rightarrow 1$ does not split.}
\item{ If $c=0$ and $\epsilon_2=-1$, then $G\cong G_b\times\bz_2\cong
D^*_{2^{n_2-1}}\times\bz_2$.}
\item{If $c=2^{n_2-1}k$ and $\epsilon_2=1$, then $G_b\cong\bz_{2^{n_2}+1}$,
and $G\cong G_b\rtimes\bz_2=\bz_{2^{n_2+1}}\rtimes\bz_2$, with $q_a^{-1}
q_bq_a=((1+2^{n_2-1})k)q$. Although $G$ contains involutions,
the extension $1\rightarrow G_a\rightarrow G
\rightarrow \lan b\ran\rightarrow 1$ does not split if $n_2>1$.}
\item {If $c=2^{n_2-1}k$ and $\epsilon_2=-1$, then $(q_aq_b)^2=1$,
so the sequence splits. }
\end{enumerate}

\end{proof}

In Section~\ref{nonlinearcase}, we will require explicit descriptions 
of the restriction maps $r^*:H^2(D_4; \bz)\rightarrow H^2(H;\bz)$
as $H$ ranges over the various subgroups of $D_4$. Rather than interrupt
the flow of that argument later, we discuss them here. The methods
used to calculate these maps are described in~\cite{Pearson}. 
Some of the specific maps are also described
there, and most of the rest were worked in the course of a
conversation with the author of that paper. 

Using the presentation $D_4\cong\lan s,t\ |\ s^4=t^2=1, tst=s^{-1}\ran$,
the subgroups (up to conjugacy) can be enumerated as follows, with
subscripts denoting generators: $G_s\approx\bz_4$, $G_{s^2}\approx
\bz_2$, $G_t\approx\bz_2$, $G_{st}\approx\bz_2$, $G_{s^2, t}\approx
\bz_2\times\bz_2$, and $G_{s^2, st}\approx \bz_2\times\bz_2$. The
integral cohomology of these groups is computed from (known)
descriptions of the cohomology with $\bz_2$ coefficients using the 
Bockstein spectral sequence. 

For a general group $K$, $H^1(K; \bz_2)
\cong\hom(K_{ ab}'\bz_2)$, and generators of $H^2(K; \bz_2)$
are often products of 1-dimensional classes. When integral cohomology
classes are lifts of powers of $\bz_2$-classes, we will name them
accordingly, even when the integral classes themselves are 
indecomposable. 

Generators of $H^1(D_4; \bz_2)\approx \hom (H_1(D_4), \bz_2)$ are given
by $e$ and $f$, where $e(\overline{s})=1$, $e(\overline{t})=0$, 
$f(\overline{s})=1$, and $f(\overline{t})=1$. This seemingly
asymmetrical choice of generators yields the convenient relation
$e\cup f=0$, while $e^2$ and $f^2$ lift to generators $\widehat{e^2}$ and
$\widehat{f^2}$ of $H^2(D_4; \bz)\approx\bz_2\times\bz_2$. 
($H_2(D_4; \bz_2)$
also contains an indecomposable element $w$ which does not lift.)

Similarly, for each copy of $\bz_2\times\bz_2$ described by an ordered
set of generators $\lan x, y\ran$, we have $H^2(\bz_2\times\bz_2; \bz)
\approx\bz_2\times\bz_2 =\lan \widehat{a^2}, \widehat{b^2}\ran$, where
$a(x)=1$, $a(y)=0$, $b(x)=0$, and $b(y)=1$. A word of warning: we use
the same notation
for $H^2$-generators of different copies of $\bz_2\times\bz_2$; meaning
can be determined from context.

Finally, $H^2(\bz_4; \bz)\approx\bz_4$, with a generator denoted by $c$.
Both  $c$ and its reduction mod 2 are indecomposable.
 
\begin{lem}\label{dfourcalcs}
The restriction maps $r^*:H^2(D_4; \bz)\rightarrow H^2(H; \bz)$,
as $H$ ranges over subgroups of $D_4$, are as follows:

\begin{center}
\begin{tabular}{r|cccccc}
&$G_{s^2, t}$&$G_{s^2, st}$&$G_{s^2}$&$G_t$&$G_{st}$&$G_s$\\
\hline
$\widehat{e^2}$& 0 &$\widehat{b^2}$& 0 & 0 & $\widehat{b^2}$&$2c$ \\
$\widehat{f^2}$&$\widehat{b^2}$& 0 & 0 & $\widehat{b^2}$ & 0 & $2c$\\
\end {tabular}\end{center}\end{lem} 

\begin{proof}
We use the commutative diagram
\begin{equation*}\begin{CD}
\lan \widehat{e^2}, \widehat{f^2}\ran\approx H^2(D_4;\bz)@>{r^*}>> H^2(H; \bz)\\
					@VVV			@VVV\\
\lan e^2, f^2, w\ran\approx H^2(D_4; \bz_2)@>{r^*}>> H^2(H; \bz_2),
\end{CD}\end{equation*}
where the vertical maps are induced by the coefficient reduction. The
right-hand map is an isomorphism for all the subgroups except for 
$\bz_4$, where $H^2(\bz_4; \bz_2)\cong H^2(\bz_4; \bz)\otimes\bz_2$.
The fact that $e^2$ and $f^2$ are squares of one-dimensional classes
makes calculation easy for every column except $G_s$. In that case,
since the generator of $H_2(\bz_4; \bz_2)$ is indecomposable, 
the map $r^*:H_2(D_4; \bz_2)\rightarrow H_2(\bz_4; \bz_2)$ is zero.
It follows immediately that each of $r^*(\widehat{e^2})$ and $r^*(\widehat{f^2})$
is either $0$ or $2c$. Determining which actually occurs requires
a closer look at the Bockstein spectral sequence. Details are 
in~\cite{Pearson}; our calculation later on actually only requires that 
each be an even multiple of $c$.
 \end{proof}

\section{The linear case}

Let us say that the group of \emph{linear} actions on $\fred$ is $W=\{A \in SO(6)\ |\ A(\fred)=\fred\}.$ What is the structure of
$W$? The homology representation $\varphi$ extends to
all of $W$. We will apply it in the proof of the following:

\begin{lem}
\label{Wstructurelemma}
$W\cong (SO(3)\times SO(3)\times \bz_2)\rtimes \bz_2$, where the 
semidirect product automorphism is a coordinate switch in the 
factors of $SO(3)\times SO(3)$.
\end{lem}

\begin{proof}
The main claim that needs to be established is the exactness of the 
following sequence:

$$1\rightarrow SO(3)\times SO(3)\stackrel{i}{\rightarrow}W
	\stackrel{\varphi}{\rightarrow}\bz_2\times\bz_2\rightarrow 1.$$

Since a rotation in $SO(3)$ is homotopic to the identity, it is clear
that any $g\in SO(3)\times SO(3)$ is in the kernel of $\varphi$.
For the reverse inclusion, suppose $\varphi(g)=e\in \bz_2\times\bz_2$.
Since $g$ acts trivially on homology, it has a fixed point
$(x,y)\in\fred$.  A fairly standard calculation shows that the
only sections of $T_{(x,y)}\fred$ with sectional curvature $K=1$
are those tangent to the two factors. Since $g$ is an isometry,
it preserves the splitting of 
$T_{(x,y)}\fred$, so it acts on $T_{(x,y)}\fred$ by a pair of 
rotations $\theta_1$ and $\theta_2$. Thus $g$ agrees with an element
of $SO(3)\times SO(3)$ on a 5-dimensional subspace of $\br^6$. Since
$\det(g)=1$,  $g\in SO(3)\times SO(3)$. This finishes the proof of 
exactness.

Now, the quotient $\bz_2\times \bz_2$ is easily seen to be generated
by $\alpha$, the isometry which acts by the antipodal map in each
factor, and $\sigma$, which switches coordinates. $\alpha$ is central
in $SO(6)$, so $\left\langle \alpha \right\rangle$ extends $SO(3)\times SO(3)$ with a
direct product. However, $\sigma$ is not central, so the product is
only semidirect. 
 
\end{proof}

For the remainder of this section, let $G$ be a finite subgroup of $W$
which acts linearly and pseudofreely and preserves orientation. 
It is well known (see~\cite{Wolf}, \cite{Kulkarni}, or 
Lemma~\ref{pairsimpliespolyhedral} below,  for example) that the finite 
groups which act pseudofreely
on $S^2$ are either cyclic, dihedral, or one of the three symmetry
groups of the Platonic solids.
These groups also act pseudofreely on $\fred$ via the
diagonal action. It will turn out that if the induced action on
homology is trivial, this is essentially, but not exactly, the only 
possibility.
If the homology action is nontrivial, things are more
complicated, and we approach them via the short exact sequence
$1\rightarrow K\rightarrow G\stackrel{\varphi}{\rightarrow} Q\rightarrow 1$.

\goodbreak
\begin{eg}
Consider the action of $\bz_5$ on $\fred$ defined by 
$\gamma\cdot(x,y)=(\gamma\cdot x, \gamma^2\cdot y),$ where $\gamma$ acts on
$S^2$ by a rotation of $2\pi /5$ around an axis. The action has four isolated
fixed points. Notice also that this action {\it resembles} the diagonal
action, in the sense that $\gamma\cdot(x,y)=(\gamma\cdot x, 
\psi(\gamma)\cdot y),$ where $\psi$ is the automorphism of $\bz_5$ 
sending $\gamma\mapsto\gamma^2.$ However, this action is {\it not}
equivalent to the diagonal action: A neighborhood of a singular point in
the quotient by the diagonal action is a cone on the lens space
$L(5,1)$, while the
corresponding neighborhood for this action is a cone on $L(5,2)$.
\end{eg}

\begin{defn}
Suppose $G$ acts on a space $X$. An action $\theta$ of $G$ on $X\times X$ 
will be called \emph{semi-diagonal} if there is
an automorphism $\psi$ of $G$ so that $\theta$ is equivalent to  
$(g, (x,y))\mapsto(gx, \psi(g)(y))$.

\end{defn}

\begin{prop}
\label{homtrivcase}
If $G$ acts trivially on homology, then $G$ is a polyhedral group,
and its action is semi-diagonal.
\end{prop}   

\begin{proof}
Since $G$
acts trivially on homology, we must have 
$G\subset SO(3)\times SO(3).$ Let $(g,h)\in G$.
Since each of $g$ and $h$ preserves orientation on $S^2$, each has a
fixed point on $S^2$. By the assumption of pseudofreeness, each fixed
set must be 0-dimensional unless $(g,h)=(e,e)$. 

{\bf Claim:} Let $\pi_1$ and $\pi_2$ denote the two obvious projections
from $SO(3)\times SO(3)\rightarrow SO(3)$. Then $\pi_1(G)\approx \pi_2(G).$

Observe that for any $g$, there is a unique $h$ so that $(g,h)\in G$.  For if
$(g,h_1)$ and $(g,h_2)$ are both in $G$ and $h_1\neq h_2$, then
$(g^{-1}, h_1^{-1})(g,h_2)=(1, h_1^{-1}h_2)\in G$. But $h_1^{-1}h_2\neq 1$,
so it has two fixed points on $S^2$. But then $\mbox{Fix}((1, h_1^{-1}h_2),
\fred)\cong S^2\times S^0$, contradicting pseudofreeness.

Now we can define $\psi:\pi_1(G)\rightarrow \pi_2(G)$ by $$g\mapsto
\mbox{ the unique $h$ such that $(g,h)\in G$.}$$
$\psi$ is clearly a homomorphism, and by symmetry, an isomorphism. 
This proves the claim. To simplify notation, write $G_1=\pi_1(G)$ and
$G_2=\pi_2(G)$.

$G_1$ acts pseudofreely on $S^2$, so it must be a polyhedral group.
Moreover, the projection $G\rightarrow G_1$ must be injective, since if
$(g,h)\in \ker\pi_1$, then $g=e$, and $h=\psi(g)=e,$ also. Similar
considerations apply to $G_2$.
Now, it is well known (see~\cite[th 2.6.5]{Wolf}) that any two 
isomorphic finite subgroups of $SO(3)$
are conjugate. That is, there is a change of 
coordinates with respect
to which we actually have $G_1=G_2$. Without loss of generality,
assume we have applied it. 
(In doing so, we have fixed, once and for all, a particular
identification between the two factors of $\fred$. For later
reference, this also defines a particular choice of coordinate
switch $\sigma:(x,y)\mapsto (y,x)$. Abstractly, $\sigma$ is only
well-defined mod $\{e\}\times SO(3)$.)
Then for any $g\in G$, we have
$g\cdot(x,y)=(gx, \psi(g)y)$. 

\end{proof}

Here it is worthwhile to note a consequence of the proof
of Lemma \ref{Platosplitslemma}: If the automorphism $\psi$
of $G$ above is inner, then after applying the coordinate
change, we may treat $\psi$ as identity.
Since $\out(\oct)$ is trivial, every $\oct$ action is diagonal.
And since $\out(\tet)\approx\out(\icos)\approx \bz_2$, each of $\tet$ and 
$\icos$ admits at most one non-diagonal action. In
fact, by embedding $\tet$ in $\oct$, we see that the outer 
automorphism of $\tet$ is realized by an $SO(3)$ conjugation,
so every $\tet$ action is diagonal. However, if there were an
$SO(3)$ conjugation which realized the outer automorphism
of $\icos$, then $S_5$ would be a subgroup of $SO(3)$. It
isn't, so there are two linear, pseudofree actions of $\icos$
which are not linearly equivalent.

If the action of $G$ on homology is nontrivial, what can we say? 
If $G\subset O(3)\times O(3)$, for example, the only
possible nontrivial action is via
$\left(\begin{smallmatrix} -1&0\\0&-1
\end{smallmatrix}\right)$. Thus we 
have a short exact sequence 
$$1\rightarrow K\rightarrow G\rightarrow \bz_2\rightarrow 0.$$
Since we know the possible groups $K$ in
the above sequence, we should expect the possible groups $G$
simply to be extensions of polyhedral groups by $\bz_2$.
But it turns out that many of
the extended actions can not be pseudofree and must be ruled out. The
following statement describes the possible groups $G$.  Explicit 
actions are constructed in its proof.

\begin{prop}\label{minusidcase}
If the homology action of $G$ is given by $\varphi(G)\approx\lan
\matrixa\ran$,
then $G\approx \bz_n\times \bz_2 $ or $G\approx\bz_{2n}$. Pseudofree
actions exist in these cases.

\end{prop}

\begin{proof}
We begin with some observations about
the geometry of $W$.

Suppose $g$ and $s$ are nontrivial elements of $SO(3)$. Each has a 
well defined axis of rotation -- we denote them respectively by
$X_g$ and $X_s$. Observe that $\mu_g(s)=s^l$ for some $l$ if and
only if $g$ leaves $X_s$ invariant. If $g$ fixes $X_s$, then $X_g=X_s$,
and $\mu_g(s)=s$. If $g$ inverts $X_s$, then $X_s\perp X_g$, and $g$ is an order 2 rotation. In this case, $\mu_g(s)=s^{-1}$. 

Lemma~\ref{Wstructurelemma} tells us that any $q\in W$ with 
$\varphi(q)=\matrixa$ has the form $q=(\alpha, \alpha)(r_1, r_2)$,
where $(r_1, r_2)\in SO(3)\times SO(3)$. Also, an element $S\in W$ 
acting homologically trivially and  pseudofreely has the form $(s_1, s_2)$, 
where neither
coordinate is $1$. It follows from the previous observations that if 
$q^2=1$ and $\mu_q(s)\ne s$ for some $s$, then $r_1$ and $r_2$ are both 
nontrivial order 2 rotations, and $q$ fixes a torus. So if $G$ acts
pseudofreely, any $q$ with $\varphi(q)=\matrixa$ and $q^2=1$ must
be central in $G$. Similar considerations show that a $q$ of order greater
than 2 leaves invariant a unique axis in each factor, so at most
one cyclic subgroup of $K$ can be normalized by $q$, and in fact, the 
generator of such a group commutes with $q$. Thus for each 
$q\in G\setminus K$, one of the following holds:
\begin{enumerate}
\item{$q$ has order 2, and $\mu_q$ is trivial.}
\item{$q$ has order greater than 2, and $q$ normalizes (at most)
a single maximal cyclic subgroup of $K$, which it centralizes.}
\end{enumerate}

Now consider the possibilities for $K$:

\noindent{\bf Case 1:} $K\cong\bz_n$.

The only groups $G$ satisfying conditions (1) and (2) are
$G\cong \bz_n\times\bz_2$ and $G\cong\bz_{2n}$. If $n$ is odd, these
two groups are isomorphic, and the extension is realized by choosing
$q=(\alpha, \alpha)\circ(1, r)$, where $r$ has order 2 and shares
an axis with the generator of $K$. If $n$ is even, this construction
realizes the case $\bz_n\times\bz_2$. $\bz_{2n}$ is realized  if
both $r_1$ and $r_2$ are order $2n$ rotations. Note that if $n$ is
odd, the latter construction still yields a $\bz_{2n}$ action, but
if $n$ is odd, then $q^n$ fixes a torus.

\noindent{\bf Case 2:} $K\cong D_n =\lan s,t|s^n=t^2=1, tst=s^{-1}\ran$.

If some $q$ exists satisfying condition (2), then $tq$ has order 2
but fails to satisfy (1). If no $q$ satisfies (2), then $n=2$, and
either $q$, $sq$, $tq$, or $stq$ has both $r_1$ and $r_2$ nontrivial,
and thus fixes a torus.

\noindent{\bf Case 3:} $K\cong\tet$, $\oct$, or $\icos$.

By Lemma~\ref{Platosplitslemma}, the sequence $1\rightarrow K\rightarrow
G\rightarrow\bz_2\rightarrow 1$ splits, so we may pick $q$ of order 2.
By condition (1), $\mu_q$ is trivial. Now let $g\in K$ be an involution.
Then $gq$ also has order 2, but $\mu_{gq}$ is nontrivial, 
contradicting (1).

\end{proof}

Next we consider the case in which $\varphi(G)\approx\lan
\matrixb\ran$. For actions which need not be locally linear,
it follows from a theorem of Bredon that:

\begin{prop}\label{bredonstheorem} 
\begin{enumerate}
\item{If the sequence $1\rightarrow K\rightarrow G\rightarrow \lan \matrixb
\ran\rightarrow 1$ splits, then $G$ can not act pseudofreely.}
\item{Suppose $G$ acts on $\fred$ with $\varphi (G)\cong\bz_2$. 
If the action of $K=\ker(\varphi)$ is pseudofree, and $1\rightarrow K
\rightarrow G\rightarrow \bz_2\rightarrow 1$ does not split, then
the action of $G$ is pseudofree.}
\end{enumerate}
\end{prop}
\begin{proof}
As a special case of Bredon's theorem~\cite[VII.7.5]{Bredon}, 
an involution $T$ on $\fred$
with $T^*\ne 1$ on $H^2(\fred ; \bz_2)$ has a fixed point set $F$
with $H^2(F; \bz_2)=\bz_2$. If the sequence  $1\rightarrow K\rightarrow 
G\rightarrow \lan \matrixb\ran\rightarrow 1$ splits, then $G$
contains such an involution. The second statement follows from the
observation that for $q\in G\setminus K$, Fix$(q^2)\subseteq$ Fix$(q)$.
 \end{proof}

\begin{cor}
If the homology action of $G$ is given by $\varphi(G)\approx\lan
\matrixb\ran$, then 
$K$ can't be $\tet$, $\oct$, or $\icos$.
\end{cor}
\begin{proof}
This follows immediately from Lemma~\ref{Platosplitslemma}.
 
\end{proof}

\begin{prop}\label{cyclicH}
If $K\cong\bz_n$  and
$\varphi(G)\approx\lan\matrixb\ran$, then either
\begin{enumerate}
\item{
	$G\cong (\bz_{m_-}\rtimes\bz_{2^{n_2+1}})\times\bz_{m_+}$, or}
\item{
	$G\cong D^*_{2^{n_2-1}m_-}\times\bz_{m_+}$,}
\end{enumerate}
with notation as in the proof of Lemma \ref{Hcycliclemma}. Pseudofree
linear actions  exist in these cases.
\end{prop}

\begin{proof}
We know from Proposition~\ref{bredonstheorem} that if $1\rightarrow
K\rightarrow G\rightarrow\bz_2\rightarrow 1$ splits, then $G$ cannot
act pseudofreely. It follows from Lemma~\ref{Hcycliclemma} that $G$
must be of type (1) or (2) (and that $n$ is even). Recall from the 
beginning of Section~\ref{grouptheory} that the structure of $G$ is
determined by the data $\gamma$ and $\gqq$. By 
Lemma~\ref{Wstructurelemma}, a hypothetical $q\in G\setminus K$ 
has the form $\sigma r$, where $\sigma :(x,y)\mapsto(y,x)$, and
$r=(r_1, r_2)\in SO(3)\times SO(3)$. Our construction of a pseudofree
$G$-action uses two ingredients:
\begin{enumerate}
\item{A choice of a particular semi-diagonal $K$-action, which
makes it possible to realize $\gamma$ as $\mu_q$.}
\item{ An appropriate choice of $r$ which ensures that $q^2=\gqq$.}
\end{enumerate}
	
As in the proof of Proposition~\ref{minusidcase}, conjugation 
by $r$ must either
fix or invert the axis of rotation of $K$ in each factor, so
$\mu_{\sigma r}=\mu_r\circ\mu_{\sigma}$ is either
$\mbox{(inversion})\circ\mu_{\sigma}$, or simply $\mu_{\sigma}$.
For simplicity, we will assume that $r_1$ and $r_2$ share the axis of
$K$, so $\mu_{\sigma r}=\mu_{\sigma}$. The construction can also be 
carried out with minor modifications if $X_{r_i}\perp X_K$. 

A semi-diagonal embedding of $K\subset SO(3)$ into $SO(3)\times SO(3)$
takes the form $k\stackrel{i}{\mapsto}(k, \psi(k))\in SO(3)\times SO(3)$
for some $\psi \in \aut (K)$. Observe that $\mu_{\sigma}(k, \psi(k))=
\sigma(k, \psi(k))\sigma=(\psi(k), k)$, and also that $i(\gamma(k))
=(\gamma(k), \psi(\gamma(k)))$. Since $\gamma^2=1$, choosing 
$\psi=\gamma$ allows us to realize $\gamma$ as $\mu_q$: 
\begin{eqnarray*}
\mu_q(i(k)) &=& \mu_{\sigma}(k, \psi(k))\\
	&=&	(\gamma(k), k)\\
	&=&	(\gamma(k), \psi(\gamma(k)))\\
	&=&	i(\gamma(k)).
\end{eqnarray*}

Next observe that $q^2=\sigma r\sigma r=(r_2r_1, r_1r_2)=
(r_1r_2, r_1r_2)$.  Thus if we choose $r_1$ and $r_2$ so that
$r_1r_2=\gqq$, then $q^2=(\gqq, \gqq)=(\gqq, \psi(\gqq))
=i(\gqq)$ (recall that $\psi=\gamma$ fixes $\gqq$).

It follows from the second part of Proposition~\ref{bredonstheorem}
that the extended action is still pseudofree.
 
\end{proof}

\begin{prop}\label{dihedralH}
If $K\cong D_n$ and $\varphi(G)\cong\lan\matrixb\ran$, then $G$
acts pseudofreely and linearly if and only if the sequence $1\rightarrow
K\rightarrow G\rightarrow\bz_2\rightarrow 1$ does not split. 
\end{prop}

\begin{proof}
We know it is necessary that the sequence not split. To prove the
converse, we assume it does not split and we construct a pseudofree
$G$-action. Let the data $\gamma$ and $\gqq$ be given. As in the
previous proof, $q=\sigma r$, and we will choose $r=(r_1, r_2)$ and 
$\psi\in \aut(D_n)$ appropriately. In this case, choose any $r_1$
and $r_2$ so that $r_2r_1=\gqq\in D_n\subset SO(3)$, and each 
$\mu_{r_i}$ leaves $D_n$ invariant. (For example, simply take 
$r_1=1, r_2=\gqq$.) Let $\psi=\mu_{r_2}\circ\gamma^{-1}\in\aut(D_n)$,
and let $i$ embed $D_n$ in $SO(3)\times SO(3)$ via $g\stackrel{i}
{\mapsto}(g, \psi(g))$. 

We show that the subgroup of $W$ generated by  $i(D_n)$ and $q$
is isomorphic to $G$, and then appeal to Proposition~\ref{bredonstheorem}
to see that the resulting action of $G$ is pseudofree. 

For the first claim, it suffices to verify that $q^2=i(\gqq)$ and
that $\mu_q\circ i=i\circ\gamma$. We have $q^2=(r_2r_1, r_1r_2)$.
On the other hand, $i(\gqq)=(r_2r_1, \psi(r_2r_1))=(r_2r_1, r_2^{-1}
(\gamma^{-1}(r_2r_1))r_2)=(r_2r_1, r_1r_2)$, because $\gamma$ fixes
$\gqq$.  Note that as a consequence, $\mu_{r_2r_1}=\mu_{\gqq}=\gamma^2$.

For $g\in D_n$, $i(\gamma(g))=(\gamma(g), \psi(\gamma(g)))=(\gamma(g),
r_2^{-1}gr_2)$. But $\mu_q(i(g))=\mu_q(g,\psi(g))=
\mu_r(\mu_{\sigma}(g,\psi(g)))=$ $\mu_r(\psi(g), g)=$ $(r_1^{-1}r_2^{-1}
\gamma^{-1}(g)r_2r_1, r_2^{-1}gr_2)=(\gamma^2(\gamma^{-1}(g)), r_2^{-1}
gr_2)=(\gamma(g), r_2^{-1}gr_2)$.

\end{proof}

\begin{prop}\label{fullcase}
Suppose $\varphi(G)=\bz_2\times\bz_2$. Then 
\begin{enumerate}
\item{$G\cong D^*_{2^{n_2}m_-}\times\bz_{m_+}$, or}
 
\item{$G\cong((\bz_{m_-}\rtimes\bz_{2^{n_2+1}}) \times\bz_{m_+})
\rtimes \bz_2, $ with $n_2>1$,}
\end{enumerate}
 with notation as in Lemma~\ref{Hcycliclemma}. Pseudofree actions 
exist in these cases.
\end{prop}
\begin{proof}
Let $a=\matrixa$ and $b=\matrixb
\in \aut (H_2(\fred))$, and suppose an extension 
$$1\rightarrow H\rightarrow G\stackrel{\varphi}{\rightarrow}\lan a \ran
\times \lan b \ran\rightarrow 1$$ exists, with $G$ acting 
pseudofreely. Let $G_a=\varphi^{-1}(\lan a\ran)$, and $G_b=
\varphi^{-1}(\lan b\ran)$.
 
By Propositions~\ref{minusidcase} and~\ref{bredonstheorem}, we know that 
$K$ must be cyclic of even order, with $G_a\cong \bz_{2n}$ or
$\bz_n\times\bz_2$, and that the sequence 
$1\rightarrow G_a\rightarrow G\rightarrow \lan b\ran\rightarrow 1$ must
not split. Lemma~\ref{technical} describes those groups $G$ which 
satisfy these criteria.

If $G_a\cong\bz_n\times\bz_2$, the geometry of $W$ puts one more 
restriction on $G$. Let $q_a^2=1$, where $\varphi(q_a)=a$. Then
$q_a=(\alpha, \alpha)(r_1, r_2)$, where $(r_1)^2=(r_2)^2=1$. If both
$r_1$ and $r_2$ are nontrivial, then $q_a$ fixes a torus. If both are
trivial, then the product of $q_a$ with the unique element of order
two in $K$ fixes a torus. So exactly one of them is nontrivial. Suppose
it is $r_1$. The element $q_b$ has the form $\sigma (s_1, s_2)$, so
$q_b^{-1}q_aq_b=(\alpha, \alpha)(1, s_2^{-1}r_1s_2)$. Thus $q_a$ and
$q_b$ can not commute in the linear case.

Lemma~\ref{technical} then shows that $G$ necessarily has one of the forms
given in the statement of the theorem. We must now construct pseudofree
actions of these groups. For simplicity, we assume all rotations
used in the constructions have the same axes.

\noindent{\bf Case 1:} $G\cong D^*_{2^{n_2}m_-}\times\bz_{m_+}$.
 
We start with an action of $G_b=D^*_{2^{n_2-1}m_-}\times\bz_{m_+}$
as provided by Proposition~\ref{cyclicH}. Extend $\psi$ (the 
semidiagonalizing automorphism of $K$) to an automorphism of $\bz_{2n}$.
Let $q_a=(\alpha, \alpha)(r, \psi(r))$, where $r^2$ is a generator 
of the Sylow 2-subgroup of $K$. Then $\lan K, q_a\ran\cong\bz_{2n}$ 
acts pseudofreely. We need only verify that $q_b^{-1}q_aq_b=q_a^{-1}$
to conclude that the subgroup of $W$ generated by $G_b$ and $q_a$ is
isomorphic to $G$: $$(s_1^{-1}, s_2^{-1})\sigma(\alpha, \alpha)(r, \psi(r))
\sigma(s_1, s_2)=(\alpha, \alpha)(\psi(r), r)=(\alpha, \alpha)(r^{-1},
\psi(r)^{-1})=q_a^{-1}.$$

\noindent{\bf Case 2:} $G\cong((\bz_{m_-}\rtimes\bz_{2^{n_2+1}}) 
\times\bz_{m_+})\rtimes \bz_2, $ with $n_2>1$.

Again, start with an action of $G_b=(\bz_{m_-}\rtimes\bz_{2^{n_2+1}})\rtimes
\bz_{m_+}$ as in Proposition~\ref{cyclicH}. Let $q_a=(\alpha, \alpha)(r, 1)$,
where $r$ is an order two rotation. Then $q_a^2=1$, so $G_a=\lan K, q_a\ran
\cong \bz_n\rtimes \bz_2$. And we have $q_b=\sigma(s_1, s_2)$, chosen
so that $(q_b)^2=k$ generates the 2-subgroup of $\bz_n$. To show that
$\lan G_b, q_a\ran\cong G$, we observe that 
$$q_a^{-1}q_bq_a=(\alpha, \alpha)(r, 1)\sigma(s_1, s_2)(\alpha, \alpha)(r, 1)
=2^{n_2-1}kq_a,$$
as required.

Pseudofreeness of the extended actions follows from
Proposition~\ref{bredonstheorem}, as usual.

\end{proof}

Gathering together the results of this section, we have:

\begin{thm}\label{biglist}
The following are all of the groups $G$ which act linearly
and pseudofreely on $\fred$:
\begin{enumerate}
\item{$\tet$, $\oct$,  $\icos$, $\bz_n$, and $D_n$, acting homologically 
trivially.}
\item{$\bz_n\times \bz_2$ and $\bz_{2n}$, where $K\approx\bz_n$, and
$\varphi(G)\approx\lan\left(\begin{smallmatrix} 
-1&0\\0&-1\end{smallmatrix}\right)\ran$.}
\item{The groups
\begin{enumerate}
\item{
	$G\cong (\bz_{m_-}\rtimes\bz_{2^{n_2+1}})\times\bz_{m_+}$, or}
\item{
	$G\cong D^*_{2^{n_2-1}m_-}\times\bz_{m_+}$,}
\item{$\lan D_n, k | k^2=\gkk, kgk^{-1}=\gamma(g)
\mbox{ for } g\in D_n\ran$, where $1\rightarrow D_n\rightarrow G
\rightarrow\bz_2\rightarrow 1$ does not split.}
\end{enumerate}
 In this case, 
$\varphi(G)\approx\lan\left(\begin{smallmatrix} 
0&1\\1&0\end{smallmatrix}\right)\ran$.}
\item{The groups
\begin{enumerate}
\item{$G\cong D^*_{2^{n_2}m_-}\times\bz_{m_+}$, or}
 
\item{$G\cong((\bz_{m_-}\rtimes\bz_{2^{n_2+1}}) \times\bz_{m_+})
\rtimes \bz_2, $ with $n_2>1$.}
\end{enumerate}
In this case, $\varphi(G)\approx\bz_2\times\bz_2$.}
\end{enumerate}

\end{thm}

This theorem can be restated in a form somewhat more amenable to 
generalization to the nonlinear case:
Let $a=\left(\begin{smallmatrix} -1&0\\0&-1\end{smallmatrix}\right)$
and $b=\left(\begin{smallmatrix} 0&1\\1&0\end{smallmatrix}\right)
\in \aut (H_2(\fred))$, and suppose an extension 
$$1\rightarrow K\rightarrow G\stackrel{\varphi}{\rightarrow}Q
\rightarrow 1$$ exists, with  $Q\subset\lan a\ran\times\lan b\ran$.
The sequence then determines the subgroups
$G_a=\varphi^{-1}(\lan a\ran)$, and $G_b=\varphi^{-1}(\lan b\ran)$. 

Conversely, suppose a homology representation of a group $G$ is described 
by a tuple $(G, K, G_a, G_b)$  as above.

\begin{cor}\label{linearclassification}
$(G, K, G_a, G_b)$ admits a linear, pseudofree action on 
$\fred$ if and only if
\begin{enumerate}
\item{ $K$ is polyhedral.}
\item{If $G_a \ne K$, then $G_a$ is abelian and $K$ is cyclic.}
\item{If $G_b\ne K$, then $1\rightarrow G_a\rightarrow G
\rightarrow \bz_2\rightarrow 1$ does not split.}
\item{If $\varphi(q_a)=a$, $\varphi(q_b)=b$, and $q_a^2=1$, then $q_a$
and $q_b$ do not commute.}
\end{enumerate}\end{cor}

\begin{proof}
 We have already seen that condition 1 is necessary and sufficient in
the homologically trivial case, and that condition 2 is necessary
and sufficient in the $\pm\left(\begin{smallmatrix}1&0\\0&1
\end{smallmatrix}\right)$ case. Condition 3 is necessary by  
Proposition~\ref{bredonstheorem} (Bredon's theorem).
On the other hand, this condition implies that
$1\rightarrow H\rightarrow G_b\rightarrow \bz_2\rightarrow 1$ does
not split, and Propositions~\ref{cyclicH} and~\ref{dihedralH} show
this to be sufficient in the $\left(\begin{smallmatrix}0&1\\1&0
\end{smallmatrix}\right)$ case. Finally, Lemma~\ref{technical} enumerates
those groups $G$ with $\varphi(G)\cong \bz_2\times\bz_2$ which
satisfy the first three conditions, and  Proposition~\ref{fullcase} shows
that the last condition is necessary and sufficient to finish the 
classification.
 
\end{proof}

\section{The nonlinear case}\label{nonlinearcase}

In this section we prove  that
any group which acts pseudofreely and locally linearly also
acts pseudofreely and linearly.
(It should be pointed out that nonlinear
actions definitely do exist: In~\cite{EdmondsEwing}, Edmonds and Ewing
construct pseudofree actions of cyclic groups with fixed point data
incompatible with linearity.)
Our strategy follows the statement of 
Corollary~\ref{linearclassification}: Bredon's theorem shows
that condition 3 is necessary. We show that conditions
1, 2, and 4 are necessary in the general case.

For the remainder of the paper, we will make constant use of the
Lefschetz fixed-point theorem and the Riemann-Hurwitz formula. We
recall them here:

Let $g: X\rightarrow X$ be a periodic, locally linear map on a
compact manifold $X$. Then $\chi(X^g)=\lambda(g)$, where 
$\lambda(g)$, the Lefschetz number of $g$,  is
the alternating sum of $g$'s traces on homology. (The formula
holds more generally, but this suffices for our purposes.) In 
particular, if $g$ acts trivially on homology, $\chi(X^g)=
\chi(X)$. In the context of homologically trivial pseudofree
actions on $\fred$, it follows that each non-identity element of
$G$ has exactly four fixed points.

The Riemann-Hurwitz formula describes the orbit structure of 
a pseudofree action of a finite group $G$ on a compact space 
$X$. If the action has singular orbits $Gx_1,\ldots, Gx_m$; 
$|G_{x_i}|=n_i$, and $|G|=N$,  then
$$\chi(X)= N\chi(X/G)- \sum_{i=1}^m (N-\frac{N}{n_i}).$$
Again, the theorem generalizes -- this time, to non-pseudofree 
actions. See, e.g. \cite{Kulkarni}. In our case ($X=\fred$; 
homologically trivial actions, for now), transfer
considerations show that $\chi(X)=\chi(X/G)=4$. Thus
$$N=\frac{4}{4-\sum_{i=1}^m (1-\frac{1}{n_i})}.$$

\begin{defn}
The {\it Riemann-Hurwitz data} for a pseudofree action is the $m$-tuple
$(n_1, \ldots, n_m)$ described in the statement of the formula; each number
$n_i$ is the size of the isotropy subgroup corresponding to one orbit.
\end{defn}

In the case of pseudofree actions of a group $G$ on $S^2$, it
is easily seen that $m\le 3$, and the only possible Riemann-Hurwitz
data are of the form $(N, N)$, $(2,2,k)$, $(2,3,3)$, $(2,3,4)$, 
or $(2,3,5)$.
By local linearity, each isotropy group is cyclic, and then group
theory calculations show that the only possible  groups are
cyclic, dihedral, tetrahedral, octahedral, or icosahedral 
(cf. \cite{Kulkarni} or Lemma~\ref{pairsimpliespolyhedral} below).
In contrast, in the case $X=\fred$, an analogous calculation
to the one in the $S^2$ case only gives the bound $m\le 7$. 
There are some 20 or so infinite families of such 
$m$-tuples satisfying the Riemann-Hurwitz formula,
 plus a finite, but quite large, number
of exceptional solutions.  Also, the tuples only describe the
sizes, and not the structures, of the isotropy groups. From this 
point of view, then, an argument directly analogous to the
one in the $S^2$ case would be intractable. However, the linear
examples exhibit two more salient features:
\begin{enumerate}
\item{ A priori, the isotropy groups might be any groups admitting
free orthogonal actions on $S^3$ (under the tangent space 
representation). In fact, they are all cyclic. }
\item{ Each element of $G$ has four fixed points: $(x,y), (x, -y),
(-x, y)$, and $(-x, -y)$.  In the cyclic,
tetrahedral, and icosahedral cases, all four fixed points always lie in
different orbits, while in the dihedral and 
octahedral cases, the fixed-point sets of some elements meet two orbits
in two points each. But all four fixed points of  $g\in G$  never 
lie in the same orbit. Thus in the list of  isotropy groups corresponding
to the Riemann-Hurwitz data, each isotropy group actually occurs
two or four times. Because of this repetition, the data in this case
correspond to the data for an action on $S^2$, where things are simpler.}
\end{enumerate}

Our strategy, then, is to show that in the general case, 
\begin{enumerate}
\item{The isotropy groups are still cyclic.}
\item{They still occur in pairs in the Riemann-Hurwitz data,}
\end{enumerate}
and then appeal to the proof in the $S^2$ case, to show:

\begin{thm}\label{groundhogdaytheorem}
A group acting pseudofreely, locally linearly, and homologically 
trivially on $\fred$ is polyhedral.
\end{thm}

The proof will use a series of lemmas. We first describe the possible
isotropy groups:

\begin{prop}\label{isotropyimpliescyclic}
Suppose a finite group $G$ acts pseudofreely, homologically
trivially, and locally linearly on a closed, simply-connected 
four-manifold $X$ with $b_2(X)\ge 2$. Then each 
isotropy group is cyclic. 
\end{prop}

\begin{proof}
It is shown in~\cite{MM2} that, without the pseudofree assumption,
such an isotropy  group $G_{x_0}$ must be abelian of rank $1$ or $2$. 
But a rank $2$ group cannot act freely on the linking sphere to $x_0$,
and thus cannot act pseudofreely on $X$.
\end{proof}

\begin{lem}\label{odivideschi}
Let $G$ and $X$ be as in Proposition~\ref{isotropyimpliescyclic}. For each
singular point $x$, let $o(X^{G_x})$ be the number of $G$-orbits
which meet $X^{G_x}$. Then $o(X^{G_x})$ divides $\chi(X)$.
\end{lem}

\begin{proof}
Since $G_x=\lan g \ran$ for some $g$, $|X^{G_x}| = |\mbox{Fix}(g)|
=\lambda(g)=\chi(X)$. 

Observe that for any $h\in G$, if $h^{-1}gh(x)=x$, then $h(x)\in 
\mbox{Fix}(g)$. In other words, $N_G(\lan g\ran)$ acts on Fix$(g)$. Two points
of Fix$(g)$ are in the same $G$-orbit if and only if they are in the
same $N_G(\lan g\ran)$-orbit, for if $gx=x$, $gy=y$, and $kx=y$, then $k^{-1}gk\in G_x =\lan g\ran$.  But each $N_G(\lan g\ran)$-orbit of the action on Fix$(g)$ has cardinality
$|\frac{N_G(g)}{\lan g\ran}|$. Since they all have the same size,
the number of orbits must divide $\chi(X)$.
 
\end{proof}

Now, corresponding to the set of singular orbits (and hence to the 
Riemann-Hurwitz data $(n_1, \ldots, n_m)$),
there is a list orbit types $G_{x_1},
\ldots, G_{x_n}$. 
We say the groups are {\it repeated in pairs} if
each orbit type occurs in this list an even number of times.

\begin{lem}\label{pairsimpliespolyhedral}
Let $G$ act homologically trivially and pseudofreely on $\fred$. If
the isotropy groups of the $G$-action are repeated in pairs, then $G$
is polyhedral.
\end{lem}

\begin{proof}
Since the groups occur in pairs, $m$ is even. Let $m'=m/2$, and 
rearrange the list so that $G_1=G_{m'+1}$, $G_2=G_{m'+2}$, and so forth. 
Since $4=N(4-\sum_{i=1}^m (1-1/n_i))$, it follows that
$2=N(2-\sum_{i=1}^{m'}(1-1/n_i))$.  As in the $S^2$ case, the 
possible $n_1, \ldots, n_{m'}$ are $(N,N), (2,2,k), (2,3,3),
(2,3,4)$, and $(2,3,5)$. In each case, these numbers represent
the sizes of maximal cyclic subgroups of $G$, and each conjugacy
class of maximal cyclic subgroups occurs in the list. 
In \cite{Kulkarni}, Kulkarni describes the groups that can correspond 
to this data. Proofs are omitted for some of his assertions, and some
details are different in our case and his, so we repeat the argument here.

$(N, N)$ clearly corresponds to a cyclic group of order $N$.

For the remaining cases, we make the following observations: If a
maximal cyclic subgroup $\lan g\ran$ of $G$ is also normal, 
it has index $\le 2$ in $G$. For $G/\lan g\ran$ operates freely
on Fix$(g)$, a set of four points. But since orbits come
in pairs, Fix$(g)$ must meet 2 or 4 $G/\lan g\ran$-orbits, so
$|G/\lan g \ran |=2$ or 1. By a similar argument, any maximal
cyclic subgroup intersecting the center of $G$ non-trivially
has index $\le 2$.

$(2,2,k)$ corresponds to a group of order $2k$. It has a cyclic,
index 2 subgroup $\lan g\ran$ of order $k$, which must be normal. 
Kulkarni assumes that his groups operate on a space with $\chi =2$, 
and uses this fact to prove that each $h\in G\setminus \lan g\ran$
must have order 2. In our case, we use the observation above.
Now assume $h\in
G\setminus \lan g\ran$ is fixed. Since each $(hg^a)$ has order 2,
$hg^ahg^a=1$, so $hg^ah^{-1}=(g^a)^{-1}$. It follows that $G$ is dihedral.

In the $(2,3,3)$ case, $|G|=12$. G is nonabelian. There are three
nonabelian groups of order 12, and two of them contain elements
of order 6. The third is $\tet$.

In the $(2,3,4)$ case, $G$ has order 24 and trivial center. With 
this in mind, a look at a table of groups (e.g. \cite[p.137]{Coxeter})
easily shows that $G=\oct$. 

Finally, in the $(2,3,5)$ case, $|G|=60$. Sylow theory shows that
G contains five copies of $D_2$, intersecting trivially, 20 
copies of $\bz_3$, and six copies of $\bz_5$. If any proper 
subgroup of $G$ were normal, it would contain every conjugate
of each element. Counting arguments show this to be impossible, 
so $G$ is simple. Thus $G=\icos$.

\end{proof}

\begin{lem}\label{repeatedinpairs}
Let $G$ act pseudofreely, locally linearly, and homologically
trivially on $\fred$. Then the isotropy groups of the $G$-action
are repeated in pairs.
\end{lem}

\begin{proof}
Say $G$ is {\em nice} if for each singular point
 $x\in \fred$, $(\fred)^{G_x}$
meets 2 or 4 $G$-orbits. If $G$ is nice, then its isotropy groups
are repeated in pairs. We assume inductively that every proper
subgroup of $G$ is nice, but that  $G$ is not. Then by 
Lemma~\ref{odivideschi}, there is some $x_0$ so that $o((\fred)^{G_{x_0}})=1$; that is, 
$G$ acts transitively on $\fred)^{G_{x_0}}$.

By Proposition~\ref{isotropyimpliescyclic}, $G_{x_0}$ is cyclic --
say $G_{x_0}=\lan g \ran$, and Fix$(g)=\{ x_0, \ldots, x_3\}$.
Let $h_1(x_0)=x_1$, $h_2(x_0)=x_2$, and $h_3(x_0)=x_3$. By 
minimality, $G$ is generated by $g, h_1, h_2$,  and $h_3$, and since
$h_i^{-1}gh_i(x_0)=x_0$, $\lan g \ran$ is normal in $G$. By minimality
again, $|g|=p$ for some prime $p$. Now, $G/\lan g \ran$ acts freely
on $\{ x_0, \ldots, x_3\}$, so $|G/\lan g \ran |=4$. Thus $G$ is an 
extension of $G_{x_0}$ by $\bz_4$ or $\bz_2\times \bz_2$. 

The remainder of the proof is an  analysis of these extensions.
Most can be ruled out by elementary group theory considerations.
The two more difficult cases use arguments essentially due to
Edmonds in~\cite{Spheres}. 

	In the following cases, consideration of the possible
automorphism actions of $H$ on $\lan g\ran$ shows that some element
of $G\setminus \lan g\ran$ must be central, and then that $g$ is
contained in a cyclic subgroup of order $2p$, contradicting minimality.
\begin{enumerate}
\item{$\bz_p\rtimes(\bz_2\times\bz_2)$, for $p>2$.}
\item{Any non-split extension of $\bz_2$ by $\bz_2\times\bz_2$.}
\item{$\bz_p\rtimes \bz_4$, for $p\equiv 3 \pmod{4}$.}
\end{enumerate}
The case $G=\bz_2\rtimes\bz_4$ must actually be abelian, so $G$ is a 
direct product.  $\bz_2\times\bz_4$
contains two cyclic subgroups of order 4. They intersect non-trivially, 
and therefore have the same fixed set. It follows that $G$ must act 
semifreely.

Two cases remain: $G=\bz_2\times\bz_2\times\bz_2$, and 
$G=\bz_p\rtimes \bz_4$, where $p\equiv 1 \pmod{4}$.

Suppose  $G=\bz_2\times\bz_2\times\bz_2$ admits a non-nice action. $G$ has seven
cyclic subgroups. Since $\bz_2\times\bz_2$ does not act freely on
$S^3$, their fixed-point sets are disjoint, and each has $\bz_2$ 
stabilizer. Since $G$ is abelian, it acts on each fixed-point set,
so each constitutes an orbit. The action has Riemann-Hurwitz data
$(2,2,2,2,2,2,2)$. 

Let $X$ be $\fred $ minus a small invariant neighborhood of the
singular set, and let $Y=X/G$. $Y$ is a compact 4-manifold with
seven $\br P^3$ boundary components $P_1, \ldots, P_7$. The 
cohomology long exact sequence for the pair $(Y, \partial Y)$
(with $\bz_2$-coefficients) shows that $\image(i^*:H^3(Y)
\rightarrow H^3(\partial Y))$ has rank 6. The covering $X\rightarrow Y$
is classified by a map $\varphi :Y\rightarrow BG$. This induces\
a map $(\varphi\circ i)^*:H^3(G)\rightarrow H^3(\bigcup_j P_j)=(\bz_2)^7$,
which factors through $H^3(Y)$. Since it factors, 
$\rk(\image(\varphi\circ i)^*)\le 6$. 

On the other hand, each $\pi_1(P_j)$ maps to a different subgroup
of $\pi_1(Y)\approx\bz_2\times\bz_2\times\bz_2$ under the natural
inclusion, and hence each $P_j$ corresponds to a different nontrivial
element of $H_1(Y)$. Since $H_1(Y)\cong\hom(H_1(Y), \bz_2)$, each 
nontrivial element of $H^1(Y, \bz_2)$ restricts non-trivially to
$H^1(P_j)$ for some $j$. By the Kunneth theorem and the cohomology
structure of $H^*(\br P^3)$, each of these has a nonzero cube
which maps to the top class of $P_j$. Thus 
$\rk(\image(\varphi\circ i)^*)=7$, a contradiction.

A somewhat similar argument covers the remaining case. Let 
$G\cong \bz_p\rtimes\bz_4$, with $p\equiv 1\pmod 4$. The  semidirect
product automorphism must have order 4; otherwise $G$ contains a
cyclic subgroup $\bz_{2p}$. Thus $G\approx\lan g,h|g^p=h^4=1,
h^{-1}gh=g^a\ran $, where  $4a=p+1$. $G$ has $p$ different 
subgroups of order 4, all of which are conjugate, so if an
action exists, it has Riemann-Hurwitz data $(4,4,4,4,p)$.  Define
$X$ and $Y$ as before. Then $Y$ has boundary consisting of five
lens spaces $L_4, L_4, L_4, L_4$, and $L_p$, with associated
inclusions $i_1, \ldots, i_5: L_{n_j}\hookrightarrow Y$. Once again,
the covering $X\rightarrow Y$ is classified by a map $\varphi :
Y\rightarrow BG$, with induced maps $\varphi\circ i_j:L_{n_j}
\rightarrow BG$. However, the cohomology calculation is just a bit 
subtler this time.

For any coefficient module $M$, the transfer map gives an isomorphism
$H^*(G;M)\rightarrow H^*(\bz_p; M)^{\bz_4}$. With $\bz_p$ coefficients,
the ring $H^*(\bz_p)$ is generated by elements $s\in H^1(\bz_p)$ and
$t\in H^2(\bz_p)$, where $t$ is the image of $s$ under the Bockstein
map. $H^*(\bz_p)$ therefore inherits a $G$-module structure from the
action of $G$ on $s$ given by $h\cdot s=as$. Thus $h\cdot t=at$, and 
$h\cdot st=(2a)st$, so the action of $h$ on $H^3(\bz_p)$ is given by 
multiplication by $-1$. This has the unfortunate consequence that
$H^3(G;\bz_p)=0$. To compensate for the $G$-action on $H^*(\bz_p)$,
we replace $\bz_p$  with a twisted coefficient module $\tp$, where 
$h$ acts by $-1$. Note that $H^*(\bz_p; \tp)\cong H^*(\bz_p;\bz_p)
\cong\bz_p$ as $\bz_p$-modules, since the restriction of the $G$-action
to its subgroup $\bz_p$ is trivial. With this twisting, we observe:
\begin{enumerate}
\item{Restriction gives an isomorphism $(\varphi\circ i_5)^*:H^3(G; \tp)
\rightarrow H^3(\bz_p; \tp)\cong \bz_p$.}
\item{For $j=1, \ldots, 4$, the maps $(\varphi\circ i_j)^*: H^3(G; \tp)
\rightarrow H^3(\bz_4; \tp)=0$ are trivial.}
\end{enumerate}

Now, since the coboundary map $\delta^*:H^3(\partial Y; \tp)\rightarrow
H^4(Y, \partial Y; \tp)$ is Poincar\'{e} dual to the augmentation
$H_0(\partial Y; \tp)\rightarrow H_0(Y; \tp)$, we see that 
$\image\ i^*:H^3(Y; \tp)\rightarrow H^3(\partial Y; \tp)$ consists
of all $(u_1, \ldots, u_5)$ in $H^3(\partial Y; \tp)$ such that
$\sum{u_j}=0$. In particular, if $u_1+u_2+u_3+u_4=0$, then $u_5=0$, 
as well. But by the observations above, there are elements $u\in
H^3(Y, \tp)$ which restrict trivially to each $H^3(L_4; \tp)$, but
non-trivially to $H^3(L_p, \tp)$. This rules out the $G$-actions in 
question.
 
\end{proof}

Theorem~\ref{groundhogdaytheorem} follows from the lemmas. Thus 
condition 1 in Corollary~\ref{linearclassification} is necessary
in the general case. In other words, any group acting pseudofreely
and homologically trivially on $\fred$ also admits a linear, 
homologically trivial, pseudofree action. 
We now proceed to prove the necessity
of condition 2.

\begin{thm}\label{matrixatheorem}
Let $G$ act pseudofreely on $\fred$, and suppose $\varphi(G)=
\lan \matrixa\ran$. Then $G$ is abelian, and $\ker (\varphi)$
is cyclic. 
\end{thm}

\begin{proof}
Let $K\subset G$ act homologically trivially, so that
$$1\rightarrow K\rightarrow G\stackrel{\varphi}{\rightarrow}
\lan\matrixa\ran \approx\bz_2\rightarrow 1$$ is exact.

\smallskip\noindent{\bf Claim:} $G/K$ acts freely on $\fred / K$. 

Let $u\in G\setminus K$, so that $\varphi(u)$ generates
$\bz_2$. If, for some $x$, $u(x)$ lies in the same orbit as $x$,
then for some $k\in K$, $ku(x)=x$. But $ku$ has Lefschetz number
$0$, so if it has a fixed point, its fixed point set must be 
2-dimensional. 

	For the same reason, $G/K$ acts freely on the set of 
singular points in $\fred / K$, and therefore it identifies the
paired orbits of Lemma~\ref{repeatedinpairs}. 
	Choose a small $G$-invariant open neighborhood $N\subset\fred$  
of the singular set of the $K$-action, and let $X=\fred\setminus N$.
Let $Y=X/G$. It follows from the  claim that $Y$ is a manifold with
boundary consisting of two or three lens spaces $L_{n_i}$. (These
$n_i$'s are exactly those which appear in pairs in the Riemann-
Hurwitz data.) 
	Since $X$ is a cover of $Y$, we can use the Cartan-Leray
spectral sequence ($E_2^{p,q}\cong H^p(G; H^q(X))
\Rightarrow H^{p+q}(Y)$) to compute $H^2(Y;\bz)$. Note that
\begin{enumerate} 
\item{Since $G$ is finite, $H^1(G; M)=0$ for any free $\bz$-module
$M$.}
\item{In general, $H^0(G;M)\cong M^G$, the submodule of $M$ fixed by $G$.
In our case, $G$ acts by $\matrixa$ on $H^2(X)$, so $H^0(G; H^2(X))=0$.}
\item{No nonzero differentials enter or leave $E^{2,0}_k\cong 
H^2(G; H^0(X))$.}
\end{enumerate}

Thus $H^2(Y)\cong H^2(G)$. Similar arguments show that, for each
component $L_{n_i}$ of $\partial Y$, $H^2(L_{n_i})\cong H^2(\bz_{n_i})$,
and that the restriction $H^2(Y)\rightarrow H^2(L_{n_i})$ is given by 
the corresponding map on subgroups. 

By Poincare duality in Y, we have an exact sequence
$$H^2(Y)\rightarrow H^2(\partial Y)\cong H_1(\partial Y)\rightarrow
H_1(Y).$$
This becomes

$$H^2(G)\stackrel{r^G_{\oplus}}{\rightarrow} \bigoplus H^2(\bz_{n_i})
\cong \bigoplus H_1(\bz_{n_i})
\stackrel{i_{\oplus}^G}{\rightarrow} H_1(G).$$ 
And since the restriction and inclusion maps factor through $K$, we 
also have a sequence
$$H^2(K)\stackrel{r^K_{\oplus}}{\rightarrow} \bigoplus H^2(\bz_{n_i})
\cong \bigoplus H_1(\bz_{n_i})
\stackrel{i_{\oplus}^K}{\rightarrow} H_1(K).$$
This sequence is not exact in general. However, from the previous
sequence, it follows that $\displaystyle{\frac{H^2(\partial Y)}
{\image (r^G_{\oplus})}}$ injects into $H_1(K)$, and that 
$\displaystyle{\frac{H^2(\partial Y)}{\image(r^K_{\oplus})}}$ 
injects into a quotient of $H_1(K)$.  The remainder of the proof 
consists of cohomology calculations showing that this is 
impossible unless $G$ is abelian and $K$ is cyclic.

\noindent{\bf Case 1:} $K$ is cyclic.

Since $G$ maps onto $\bz_2$, $[G, G]\subseteq K$. Thus if $G$ is 
nonabelian, the kernel of $\lan k\ran\rightarrow H_1(G)$ is nontrivial.
In the exact sequence
$$H^2(G)\stackrel{r^G_{\oplus}}{\rightarrow}  H^2(\bz_n)\oplus
H^2(\bz_{n})\cong  H_1(\bz_n)\oplus H_1(\bz_n)
\stackrel{i_{\oplus}^G}{\rightarrow} H_1(G),$$
the image of $r^G_{\oplus}$ lies in the diagonal subgroup of
$H^2(\bz_n)\oplus H^2(\bz_n)$. If $G$ is nonabelian, the kernel
of the inclusion map does not.
\smallskip

\noindent{\bf Case 2:} $K=\tet$, $\oct$, or $\icos$.

To every polyhedral group $K$, there corresponds a binary polyhedral
group $\widetilde{K}$ such that $1\rightarrow \bz_2\rightarrow 
\widetilde{K}\rightarrow K\rightarrow 1$ is exact. The 
Lyndon-Hochschild-Serre spectral sequence ($E^{p,q}_2\cong 
H^p(K; H^q(\bz_2))\Rightarrow H^{p+q}(\widetilde{K})$, in 
this case) relates the cohomologies of the groups in this 
sequence. It shows, in particular, that $H^2(K)$ injects into
$H^2(\widetilde {K})$. But since each $\widetilde{K}$ acts freely 
on $S^3$, Poincare duality shows that $H^2(\widetilde{K})\approx
H_1(\widetilde{K})$. Thus $H^2(K)\hookrightarrow H_1(\widetilde{K})$.

By computing the sizes of $H_1(\widetilde{K})$ in each case, 
we find that $|H^2(\tet)|\le 3$, $|H^2(\oct)|\le 2$, and 
$H^2(\icos)=0$. On the other hand, $|\bigoplus H^2(\bz_{n_i})|=
n_1 \times n_2 \times n_3$. Finally, $H_1(\tet)=\bz_3$, 
$H_1(\oct)=\bz_2$, and $H_1(\icos)=0$. In each case, 
$\displaystyle{\frac{H^2(\partial Y)}{\image (r^K_{\oplus})}}$ 
is too large to inject into a quotient of $H_1(K)$. 
\smallskip

\noindent{\bf Case 3:} $K=D_n$. 

Recall that

\begin{equation*}
H_1(D_n)\approx
	\begin{cases}
	\bz_2=\lan\overline{t}\ran \text{if $n$ is odd,}\\
	\bz_2\times\bz_2 =\lan \overline{s}, \overline{t}\ran
		\text{if $n$ is even,}
	\end{cases}
\end{equation*}
and 
\begin{equation*}
H^2(D_n)\approx 
	\begin{cases}
	\bz_2 \text{ if $n$ is odd,}\\
	\bz_2\times\bz_2 \text{ if $n$ is even}.
	\end{cases}
\end{equation*}
Since $n_1=n_2=2$, and $n_3=n$, $|\bigoplus H^2(\bz_{n_i})|=4n$.
Except in the cases $n=2$ and $n=4$, it is immediate that 
$\displaystyle{\frac{H^2(\partial Y)}{\image (r^K_{\oplus})}}$
cannot inject into any quotient of $H_1(D_n)$. We consider
the remaining possibilities:

Suppose $n=2$. Using Lemma~\ref{dihedralextensionlemma}, we find that
the only nonabelian extension $1\rightarrow D_2\rightarrow G\rightarrow
\bz_2\rightarrow 1$ has 
$G= D_4\cong\lan s,t|s^4=t^2=1, tst^{-1}=s^{-1}\ran$, with 
$D_2=\lan s^2, t\ran$. Thus we have an exact sequence
 
$$H^2(D_4)\stackrel{r^{D_4}_{\oplus}}{\rightarrow}  H^2(\bz_2)\oplus
H^2(\bz_2)\oplus H^2(\bz_2)
\cong  H_1(\bz_2)\oplus H_1(\bz_2)\oplus H_1(\bz_2)
\stackrel{i_{\oplus}^{D_4}}{\rightarrow} H_1(D_4),$$
where the three $\bz_2$ subgroups are generated by $s^2$, $t$, and
$s^2t$. Since the restrictions and inclusions factor through $D_2$,
Lemma~\ref{dfourcalcs} shows that $|\image (r^{D_4}_{\oplus})|=
|\ker (i^{D_4}_{\oplus})|=2$. This contradicts exactness. 

	$D_2$ also has two abelian extensions: $1\rightarrow D_2
\rightarrow D_2\times\bz_2\rightarrow\bz_2\rightarrow 1$, and
$1\rightarrow D_2\rightarrow \bz_4\times\bz_2\rightarrow\bz_2
\rightarrow 1$. In each case, the restriction and inclusion maps can
be explicitly calculated, and the sequence is seen not to be
exact. $D_2$ is the only  candidate for 
$\ker (\varphi)$ which is abelian, but not cyclic.

Finally, suppose $K=D_4$. Then we have $\bz_{n_1}=\lan t\ran\approx
\bz_2$, $\bz_{n_2}=\lan st \ran\approx\bz_2$, and $\bz_{n_3}=\lan s \ran
\approx \bz_4$. Consulting Lemma~\ref{dfourcalcs} again, we see that 
$r^{D_4}_\oplus$ has matrix $\left(\begin{smallmatrix} 0&1\\1&0\\2&2
\end{smallmatrix}\right)$ relative to the bases discussed there.
And $i^{D_4}_{\oplus}$ has matrix $\left(\begin{smallmatrix} 0&1&1\\
1&1&0\end{smallmatrix}\right)$, as is easily checked.

Now, if $r^G_{D_4}$ is not onto, then the element counts which 
applied for most $n$ also apply for $D_4$. But if $r^G_{D_4}$ is
onto, then $\displaystyle{\frac{\oplus H^2(\bz_{n_i})}
{\image(r^{D_4}_{\oplus})}}$ should inject into $H_1(D_4)$, and hence
we should have $\ker(i^{D_4}_{\oplus})\subseteq\image(r^{D_4}_{\oplus})$.
However, the element $(1,1,1)\in H_1(\lan t\ran)\oplus H_1(\lan st\ran)
\oplus H_1(\lan s\ran)$ is in the kernel, but not in the image.

\end{proof}

Condition 3 of Corollary~\ref{linearclassification} is necessary by
Bredon's theorem. A related result of Bredon helps us establish
the necessity of condition 4, and thus complete the proof of 
the main theorem:

\begin{prop}\label{lastcondition}
Suppose $G$ acts pseudofreely, $\varphi(q_a)=a$, $\varphi(q_b)=b$,
and $q_a^2=1$. Then $q_a$ and $q_b$ do not commute.
\end{prop}
\begin{proof} 

Given the necessity of conditions 1, 2, and 3 of 
Corollary~\ref{linearclassification}, it suffices to rule out
the possibility that $G=\lan q_a, q_b\ran\cong \bz_2\times \bz_{4k}$, where 
$k\ge 1$.  In the case of the linear models (Propostion~\ref{fullcase}),
the argument divided into two parts: If $q_a=(\alpha, \alpha)$, then 
$q_aq_b^{2k}$ fixes a torus, contradicting pseudofreeness. And if
$q_a =(\alpha,\rho)$ or $(\rho, \alpha)$, then the fact 
that $q_b$ exchanges factors of $\fred$
means that $q_a$ and $q_b$ cannot commute, so no action exists, even
with a two-dimensional singular set. 

	Assume, then, that $q_a$ and $q_b$ commute.  Now, $q_a$ must act
freely, so $X=\fred /\lan q_a\ran$ is a manifold which inherits an action
by $\lan q_b\ran=\bz_{4k}$. 
As motivation, we note that the linear models for $q_a$ are distinguished
by the intersection forms (with $\bz_2$ coefficients) of the quotient
spaces:
In $X_1=\fred/(\alpha, \alpha)$, the diagonal $S^2$ maps to an embedded
$\br P^2$ with self-intersection $1$, so $w_2(X_1)\ne 0$. On the other
hand, generators of second homology for $X_2=\fred/ (\alpha, \rho)$
are given by  the images of $S^2\times *$, 
(where $*$ is some point fixed by $\rho$), and 
$(*\ \cup\ -*)\times S^2$, (where $*$ is any point in the first factor). 
Each of these has trivial self-intersection, so $w_2(X_2)=0$. 
To complete the proof, we will show (just as in the linear
models) that if $\lan q_b\ran$ acts on $X$,
then $w_2(X)\ne 0$. Bredon's theorem then guarantees a two-dimensional
singular set.

Let $x$ and $y$ denote the standard generators for $H^2(\fred;\bz_2)$, 
and consider the cohomology spectral sequence of the covering:
$E_2^{p,q}(X)=H^p(\bz_2; H^q(\fred;\bz_2))\Rightarrow H^{p+q}(X;\bz_2)$.

\begin{table}[h]
\caption{$E_2(X)$}\label{L=0table}
$\begin{array}{c|lllllll}
\lan xy\ran & \bz_2 &\bz_2 &\bz_2 &\bz_2 &\bz_2 &\bz_2&\cdots \\
 & 0&0&0&0&0&0 \\
\lan x,y\ran & \bz_2\oplus\bz_2&  \bz_2\oplus\bz_2&  \bz_2\oplus\bz_2&  \bz_2\oplus\bz_2&  \bz_2\oplus\bz_2&  \bz_2\oplus\bz_2&\cdots \\
&0&0&0&0&0&0\\
\lan 1\ran & \bz_2 &\bz_2 &\bz_2 &\bz_2 &\bz_2 &\bz_2&\cdots \\
\hline
 & \bz_2 & \lan s\ran & \lan s^2\ran & \lan s^3\ran & \lan s^4\ran& 
\lan s^5\ran& \cdots\\ 
\end{array}$\end{table} 

It follows from the multiplicative structure of the spectral sequence that
the behavior of the entire
$E_2$ page is determined by $d_2(x)$ and $d_2(y)$.
At least one must be nonzero, since the sequence converges to the cohomology
of a four-manifold. And since $\lan q_b\ran$ acts on the quotient,
the differentials must respect the induced action of $q_b$ on $H^2(\fred)$.
Hence $d_2(x)=d_2(y)=1s^3$, and $\ker d_2$ is generated by $x+y$. It is
easy to check that $E_3=\ldots=E_{\infty}$, so $u=[x+y]$ is a permanent
cocycle.

	This spectral sequence is identified by a homotopy equivalence
with that of the fibration $\fred\to E\bz_2\times_{\bz_2}(\fred)
\to B\bz_2$ (cf.~\cite[IX.15]{Hu}), and under this identification,
the cocycles which live to $E_{\infty}^{0,*}$ 
are those in the image of $p^*:H^*(X)\to H^*(\fred)$.
Thus there is a class $u\in H^2(X;\bz_2)$ which lifts to 
$x+y\in H^2(\fred;\bz_2)$. 
 
We claim that $u\cup u\ne 0$. To see this from the homological point
of view, let $C_*(X)$ denote the singular  chain complex
of $X$. Then $\fred =\widetilde{X}$, and $C^*(\widetilde{X};\bz)\cong
C^*(X; \bz)\otimes_{\bz}\bz[\bz_2]$. 
The covering projection induces $p^*: C^*(X)\to
C^*(\widetilde{X})$ via $c\mapsto c\otimes(1+q_a)$. 
Moreover, because $C_*(X)$ is a free $\bz$-module, there is a natural
isomorphism $\mu:C^*(X;\bz)\otimes\bz_2\to C^*(X; \bz_2)$, so that
$H_*(C^*(X)\otimes\bz_2)= H^*(X; \bz_2)$.

Every class in $H^2(\fred;\bz_2)$ is integral. Thus we can choose
a cochain $\upsilon\in C^*(X; \bz)$ so that $\delta(p^*(\upsilon))=
\delta(\upsilon\otimes(1+q_a))= 0$
in $C^3(\tilde{X};\bz)$ and such that
$\upsilon\otimes 1$ is a cocycle representing
$u$ in $H^2(X;\bz_2)$. (Note that $\upsilon$ itself need not be an integral 
cocycle. For related discussion, see~\cite{AcostaLawson}.)
Then $[p^*(\upsilon)\otimes 1] = x+y\in H^2(\widetilde{X};\bz_2)$, and
$[p^*(\upsilon)]\in H^2(\widetilde{X};\bz)$ represents an integral lift
of $x+y$ -- say $(2m+1)\hat{x}+(2n+1)\hat{y}$. So $[p^*(\upsilon)\cup
p^*(\upsilon)] = (8mn+4m+4n+2)(\hat{x}\cup\hat{y})\equiv 2\mod{4}$.

But $p^*(\upsilon)\cup p^*(\upsilon)$ is also an equivariant cochain
and hence must be of the form $\gamma\otimes(1+q_a)$, where $\gamma\in
C^4(X;\bz)$. So $[\gamma]\equiv 1\mod 2$, and $[\gamma\otimes 1]$ represents
$ u\cup u$. Thus $u\cup u\ne 0$, so $w_2(X)\ne 0$.

A heuristic geometrical argument is more direct: the Poincar\'e dual 
of $u$ can be represented by a surface $F$, and 
$\widetilde{F}\subset \fred$ will represent
an integral lift of $PD(x+y)$, so $\widetilde{F}\cdot\widetilde{F}
\equiv 2\pmod 4$. The intersection points will be paired up by
the $\lan q_a \ran$-action, so $F\cdot F\equiv 1\pmod 2$, 
and $u\cup u\ne 0$.



\end{proof}

Combining Corollary~\ref{linearclassification}, 
Theorem~\ref{groundhogdaytheorem}, Theorem~\ref{matrixatheorem},
Proposition~\ref{lastcondition},
and Proposition~\ref{bredonstheorem}, we have our main theorem:
\begin{thm}\label{maintheorem}
Any finite group  which admits a locally linear, orientation 
preserving, pseudofree action on $\fred$ also admits a linear,
orientation preserving, pseudofree action.
\end{thm}
 
In fact, the pairs (group, cohomology representation) are exactly those
which occur in the linear case.

\bibliographystyle{abbrv}
\bibliography{mybiblio}

Ê
\end{document}